\documentclass[reqno]{amsart}

\usepackage{setspace}

\usepackage{hyperref,todonotes,comment}
\usepackage[color,matrix,arrow]{xy}
\usepackage{graphics,graphicx,amssymb}
\usepackage{pagecolor,lipsum}
\usepackage{latexsym}
\usepackage{mathrsfs}
\xyoption{curve}

\definecolor{color}{HTML}{8E1900}
\usepackage{hyperref}
\hypersetup{ 
colorlinks,
citecolor=color,
filecolor=color,
linkcolor=color,
urlcolor=color
}

\newtheorem{lemma}{Lemma}[section]
\newtheorem{proposition}[lemma]{Proposition}
\newtheorem{theorem}[lemma]{Theorem}
\newtheorem{corollary}[lemma]{Corollary}

\newtheorem*{theoremA}{Theorem}

\theoremstyle{definition}

\newtheorem{definition}[lemma]{Definition}
\newtheorem{remark}[lemma]{Remark}

\newcommand{\mfk}[1]{\mathfrak{#1}}

\newcommand{\msf}[1]{\mathsf{#1}}
\newcommand{\msc}[1]{\mathscr{#1}}
\newcommand{\mbf}[1]{\mathbf{#1}}
\newcommand{\opn}[1]{\operatorname{#1}}
\newcommand{\ot}{\otimes}

\DeclareMathOperator{\Hom}{Hom}
\DeclareMathOperator{\End}{End}

\DeclareMathOperator{\Ext}{Ext}

\DeclareMathOperator{\rep}{rep}
\DeclareMathOperator{\Rep}{Rep}

\DeclareMathOperator{\res}{res}
\DeclareMathOperator{\SL}{SL}
\DeclareMathOperator{\Rind}{R.\!ind}
\DeclareMathOperator{\Lind}{L.\!ind}

\newcommand{\dG}{\check{G}}

\renewcommand{\1}{\mathbf{1}}
\renewcommand{\O}{\mathscr{O}}

\renewcommand{\hat}{\widehat}
\renewcommand{\binom}[2]{{\Small\left[\begin{matrix}\ #1\ \\ #2 \end{matrix} \right]}}

\definecolor{page_color}{HTML}{000000}
\definecolor{text_color}{HTML}{F0EAD6}

\makeindex

\title[]{Revisiting the Steinberg representation at arbitrary roots of $1$}
\date{\today}

\author{Cris Negron}
\address{Department of Mathematics, University of Southern California, Los Angeles, CA 90007}
\email{cnegron@usc.edu}

\begin{document}

\maketitle

\begin{abstract}
We consider quantum group representations for a semisimple algebraic group $G$ at a complex root of unity $q$.  Here $q$ is allowed to be of any order.  We revisit some fundamental results of Parshall-Wang and Andersen-Polo-Wen from the 90's.  In particular, we show that the category $\Rep(G_q)$ of quantum group representations has enough projectives and injectives, and that a $G_q$-representation is projective (resp.\ injective) if and only if its restriction to the small quantum group is projective (resp.\ injective).  Our results reduce to an analysis of the Steinberg representation in the simply-connected setting, and are well-known at odd order $q$ via works of the aforementioned authors.  The details at arbitrary $q$ have, to our knowledge, not appeared in the literature up to this point.
\end{abstract}

\setcounter{tocdepth}{1}

\section{Introduction}

Let $k$ be an algebraically closed field of characteristic $0$.  For a semisimple algebraic group $G$ and a choice of nonzero parameter $q$ in $k$, we have the associated category of quantum group representations $\Rep(G_q)$.  Equivalently, we consider the category of integrable, appropriately graded representations for Lusztig's divided power quantum enveloping algebra \cite{lusztig90II}.  We are especially interested in the case where $q$ is a root of unity.
\par

The main point of this paper is to recover some well-known phenomena for quantum groups, but now at arbitrary quantum parameters.  Specifically, we want to study projective and injective representations in $\Rep(G_q)$, restrictions of such representations to the small quantum group, and the behaviors of induction from the small quantum group to the big quantum group.  Our analyses hinge on a study of the Steinberg representation in $\Rep(G_q)$ for simply-connected $G$.
\par

To place things in their proper historical context; in the late 80's and early 90's there were a number of papers which established the fundamental behaviors of quantum group representations at \emph{odd} order quantum parameters, or parameters of odd prime power order.  One can see, for example, works of Andersen-Polo-Wen, Parshall-Wang, and Lusztig from this time period \cite{lusztig89,parshallwang91,andersenpolowen91,andersenpolowen92,andersen92}.
\par

In the current mathematical and physical landscape however, we find ourselves in need of a deeper understanding of quantum group representations, specifically at \emph{even} order parameters.  This is due to recent developments in logarithmic conformal field theory and, what one might call, non-compact topological field theory.  From a mathematical perspective, the logarithmic, or non-compact modifier indicates a non-semisimplicity in the theory.  For appearances of quantum groups in such conformal and topological settings one can see the essentially-random samplings \cite{feiginetal06,pearceetal06,gainutdinovetal06,bfgt09,fuchsschweigert17,feigintipunin,sugimoto21,gannonnegron} and \cite{costellocreutziggaitto19,creutzigdimoftegarnergeer,gukovetal21,schweigertwoike21,derenzietal,costantinogukovputrov}, respectively.  For a specific point of reference, one can consider the type $A$ claims from the work \cite{creutzigdimoftegarnergeer} and the implications such claims have for quantum group representations at a $2p$-th root of $1$ (see also \cite{kapustinwitten07,gaiotto19,frenkelgaiotto20}).
\par

In the present paper we recover many of the results from the foundational texts \cite{lusztig89,parshallwang91,andersenpolowen91,andersenpolowen92,andersen92}, now at completely general $q$ and completely general $G$.  To begin, we have the following.

\begin{theoremA}[{\ref{thm:steinberg}, \ref{thm:steinberg_semis}}]
Consider a semisimple algebraic group $G$ at an arbitrary root of unity $q$.
\begin{enumerate}
\item The category $\Rep(G_q)$ has enough projectives and injectives, as does its subcategory $\rep(G_q)$ of finite-dimensional representations.
\item An object in $\rep(G_q)$ is projective (resp.\ injective) if and only if it is projective (resp.\ injective) in $\Rep(G_q)$.
\item $\Rep(G_q)$ is Frobenius, in the sense that its projectives and injectives agree.
\end{enumerate}
Furthermore, when $G$ is simply-connected, the Steinberg representation $\opn{St}=L(\rho_l)$ is projective and injective in $\Rep(G_q)$.
\end{theoremA}

Here the Steinberg representation $\opn{St}$ is the distinguished simple representation of highest weight $\rho_l=\frac{1}{2}\sum_{\gamma\in \Phi^+}(l_\gamma-1)\gamma$, where $l_\gamma=\opn{ord}q^{(\gamma,\gamma)}$.  This representation has proved fantastically useful in studies of both quantum and modular representation theory, and our work is no exception.  Indeed, one establishes points (1)--(3) above by first considering the behaviors of the Steinberg representation in the simply-connected setting.

We note that the results of Theorems \ref{thm:steinberg} and \ref{thm:steinberg_semis} first appeared in works of Parshall-Wang and Andersen-Polo-Wen \cite{parshallwang91,andersenpolowen91}, at odd order $q$, and are deducible by experts at arbitrary $q$ \cite{andersen03,andersen}.
\par

Now, at this point, we would like to make a seemingly benign assertion.  Namely, we propose that the Steinberg representation restricts to a simultaneously projective, injective, and simple representation over the small quantum group.  However, both for this statement and for the mathematics which precedes it, one first has to decide what the small quantum group actually is, or should be, at arbitrary $q$.
\par

Addressing this point is another implicit goal of this work, as well as its predecessor \cite{negron21} and successor \cite{negron}.
\par

The answer to this question is a little complicated, as there doesn't seem to be one small quantum group, but a family of small quantum groups which are parametrized by central subgroups in $G$, or more precisely by central subgroups in Lusztig's dual group $G^\ast$ from \cite{lusztig93}.  Here we work with the smallest version of the small quantum group, which we refer to as the ``smallest quantum algebra''\footnote{Most of our results for $\bar{u}_q$ will hold for any reasonable version of the small quantum group, at arbitrary $q$.  See Section \ref{sect:compare} for some discussion.}
\[
\bar{u}_q=\bar{u}(G_q).
\]
\par

At odd order $q$ and simply-connected $G$ the algebra $\bar{u}_q$ is just Lusztig's original finite-dimensional Hopf algebra from \cite{lusztig90,lusztig90II}.  At even order $q$, $\bar{u}_q$ is some variant of Arkhipov and Gaitsgory's small quantum group from \cite{arkhipovgaitsgory03}, and across all $q$ one recovers the representation category $\Rep(\bar{u}_q)$ uniformly as the fiber of $\Rep(G_q)$ along quantum Frobenius.  (See Section \ref{sect:compare} for more details on this point.)
\par

In terms of its formal structures, there is a surjective map from the quantum function algebra to the dual $\O(G_q)\to \bar{u}^\ast_q$ which makes $\bar{u}^\ast_q$ into a right $\O(G_q)$-module coalgebra.  Equivalently, at the level of categories, $\Rep(\bar{u}_q)$ has the structure of a module category over $\Rep(G_q)$ and the restriction map
\[
\res^-_q:\Rep(G_q)\to \Rep(\bar{u}_q)
\]
has the structure of a $\Rep(G_q)$-module category functor.  The kernel of this functor is precisely the image of Lusztig's quantum Frobenius $\opn{Fr}:\Rep(G^\ast_\varepsilon)\to \Rep(G_q)$.
\par

Descriptions of $\bar{u}_q$ in terms of its generators, dimension, triangular decomposition, braid group action, etc., can be found in Sections \ref{sect:sqa} and \ref{sect:sqa_at_lacing} below.  Relations with quantum Frobenius are discussed in Section \ref{sect:normal}.
\par

Our choice of $\bar{u}_q$, as the ``cokernel" of Lusztig's quantum Frobenius functor, facilitates a Steinberg decomposition for simple representations over $G_q$.

\begin{theoremA}[\ref{thm:lus_st}/\ref{cor:lus_st}]
Suppose $G$ is simply-connected.  Then every simple representation $L(\lambda)$ in $\Rep(G_q)$ admits a unique decomposition
\[
L(\lambda)=L(\lambda')\ot L(\lambda''),
\]
where $L(\lambda')$ is a simple representation which is in the image of quantum Frobenius and $L(\lambda'')$ is a simple representation which has simple restriction to $\Rep(\bar{u}_q)$.
\end{theoremA}

This result is known at odd order parameters through work of Lusztig \cite{lusztig89}, and has its origins in modular representation theory \cite{steinberg63}.  We also have the following, which generalizes results from Andersen-Polo-Wen \cite{andersenpolowen92}.

\begin{theoremA}[{\ref{thm:steinberg}/\ref{thm:steinberg_semis}}]
An object in $\Rep(G_q)$ is projective (equivalently injective) if and only if its restriction to $\Rep(\bar{u}_q)$ is projective (equivalently injective).
\end{theoremA}

In establishing the findings of Theorems \ref{thm:steinberg} and \ref{thm:steinberg_semis} we also show that the induction functor from $\Rep(\bar{u}_q)$ to $\Rep(G_q)$ is exact and faithful.  This result appears in Theorem \ref{thm:ind_ex} below.
\par

For a more elaborate historical accounting of the findings herein, and their many predecessors in the literature, we invite the reader to see Section \ref{sect:history_inc} at the conclusion of the text.  We conclude the introduction with a short discussion of the smallest quantum algebra and its relation to ``the" small quantum group, as constructed in the works \cite{creutziggainutdinovrunkel20,gainutdinovlentnerohrmann,negron21,gannonnegron,negron} for example.

\subsection{The smallest quantum algebra vs.\ the small quantum group}
\label{sect:compare}

The small quantum group for $G$ at $q$, whatever it might be \cite{arkhipovgaitsgory03,creutziggainutdinovrunkel20,gainutdinovlentnerohrmann,gaitsgory21,negron21,negron}, is represented by an algebra $u_q$ which fits into an exact sequence
\[
k\to \bar{u}_q\to u_q\to k\Sigma\to k.
\]
Here $u_q$ is a flat extension of $\bar{u}_q$, $\Sigma$ is a finite abelian group, and the final map induces an isomorphism $k\ot_{\bar{u}_q}u_q=u_q\ot_{\bar{u}_q}k\overset{\sim}\to k\Sigma$.  From this sequence one sees that the representation theory, and the homological algebra, for $u_q$ are generically identified with that of $\bar{u}_q$.  For example, a $u_q$-representation is injective if and only if it has injective restriction to $\bar{u}_q$.  In this way the conclusions of Theorem \ref{thm:steinberg_semis} transfer immediately from $\bar{u}_q$ to $u_q$.  However, highly refined structural results, such as the Steinberg decomposition from \cite{lusztig89} or Theorem \ref{thm:lus_st}, only hold for $\bar{u}_q$.
\par

From a categorical perspective, the small quantum group represents the fiber
\begin{equation}\label{eq:182}
\opn{Vect}\ot_{\Rep(\dG)}\Rep(G_q)\overset{\sim}\longrightarrow \Rep(u_q)
\end{equation}
of $\Rep(G_q)$ along a central map $\mbf{Fr}:\Rep(\dG)\to \Rep(G_q)$ which identifies the Tannakian center in $\Rep(G_q)$ with representations of a dual group $\dG$ \cite{arkhipovgaitsgory03,gaitsgory21,negron}.  This fiber construction annihilates the Tannakian center in $\Rep(G_q)$, and realizes $\Rep(u_q)$ as the unique (super-)modularization of the braided category of $G_q$-representations.  This makes $u_q$ a natural object to study from the perspective of quantum algebra or quantum topology.
\par

On the other hand, representations for the smallest quantum algebra are identified with the fiber along Luztig's usual quantum Frobenius
\begin{equation}\label{eq:181}
\opn{Vect}\ot_{\Rep(G^\ast_{\varepsilon})}\Rep(G_q)\overset{\sim}\longrightarrow \Rep(\bar{u}_q)
\end{equation}
\cite{arkhipovgaitsgory03,negron}.  Since $\Rep(G^\ast_\varepsilon)$ is generally non-central in $\Rep(G_q)$, and often not even symmetric, this fiber is most naturally understood as a right module category over $\Rep(G_q)$.\footnote{We could have taken the fiber on the right, and had an action of $\Rep(G_q)$ on the left, but it doesn't matter.  These two constructions are equivalent.}
\par

It is at the level \eqref{eq:181} that we recover classical structure theorems which ``factor" the representation theory of $\Rep(G_q)$ between that of $\Rep(G^\ast_{\varepsilon})$ and $\Rep(\bar{u}_q)$, as in \cite{lusztig89,parshallwang91,andersenpolowen91,andersenpolowen92}.

\subsection{Acknowledgments}

Thanks to Henning Andersen, Gurbir Dhillon, Pavel Etingof, Eric Friedlander, Dennis Gaitsgory, Simon Lentner, and Julia Pevtsova for helpful discussions on many topics which are related to this work.  Thanks also to H.\ Andersen and S.\ Lentner for comments on an earlier version of this text.  The author was supported by NSF Grant No.\ DMS-2149817 and Simons Collaboration Grant No.\ 999367.

\tableofcontents

\section{Preliminaries on categories and Hopf algebras}

We fix an algebraically closed field $k$ of characteristic $0$. All categories are $k$-linear categories, algebras are $k$-algebras, coalgebras are $k$-coalgebras, etc.  We let $\opn{Vect}$ denote the tensor category of arbitrary $k$-vector spaces.

\subsection{Representations and corepresentations}

Let $B$ be an algebra.  By a $B$-module we always mean a left $B$-module.  We let $\Rep(B)$ denote the category of all $B$-modules $V$ which are the unions of their finite-dimensional submodules.  Dually, for any coalgebra $C$, by a $C$-comodule we always mean a right $C$-comodule.  We let $\opn{Corep}(C)$ denote the category of arbitrary $C$-comodules, and note that any comodule is already the union of its finite-dimensional subcomodules (see for example \cite{montgomery93}).

\subsection{Tensor categories and module categories}

A tensor category is an abelian, compactly generated monoidal category $\msc{A}$ whose compact and rigid objects agree.  We require additionally that the subcategory $\msc{A}^c$ of compact/rigid objects in $\msc{A}$ is an essentially small, locally finite abelian subcategory \cite[Definition 1.8.1]{egno15}, that the unit object in $\msc{A}$ is simple, and that the tensor product commutes with colimits in each factor.  These generation restrictions imply, for example, that any tensor category $\msc{A}$ is a Grothendieck abelian category \cite[Theorem 8.6.5]{kashiwaraschapira05}.  Since tensoring with any rigid object is exact \cite[proof of Proposition 4.2.1]{egno15}, and filtered colimits are exact in any Grothendieck category, it follows that the product $\ot$ on any tensor category is exact in each factor.
\par

The following lemma allows us to do homological algebra in any tensor category $\msc{A}$, in the expected ways.

\begin{lemma}\label{lem:inj_A}
Any tensor category $\msc{A}$ has enough injectives.  Furthermore, for any injective $I$                 and arbitrary $V$ in $\msc{A}$, the products $I\ot V$ and $V\ot I$ are also injective.
\end{lemma}

\begin{proof}
Any Grothendieck abelian category has enough injectives \cite[Theorem 9.6.2]{kashiwaraschapira05}, so that the first statement follows.  Now, for rigid $V$ and arbitrary injective $I$, we have natural identifications
\[
\Hom_{\msc{A}}(-,V\ot I)=\Hom_{\msc{A}}(V^\ast\ot-,I)\ \ \text{and}\ \ \Hom_{\msc{A}}(-,I\ot V)=\Hom_{\msc{A}}(-\ot {^\ast V},I).
\]
Via exactness of the product functor $\ot$, it follows that $V\ot I$ and $I\ot V$ are injective.  Finally, since $\msc{A}$ is generated by compact objects, a version of Baer's criterion implies that filtered colimits of injectives in $\msc{A}$ remain injective.  This implies that $V\ot I$ and $I\ot V$ are injective at arbitrary $V$.
\end{proof}

A (right) module category over a tensor category $\msc{A}$ is an abelian, compactly generated category $\msc{M}$ which is equipped with an associative action of $\msc{A}$ on the right of $\msc{M}$, in the precise sense of \cite[Definition 7.1.1]{egno15}.  We require that the subcategory of compact objects $\msc{M}^c$ is an essentially small, locally finite, abelian subcategory in $\msc{M}$ and that the action map $\ot:\msc{M}\times\msc{A}\to \msc{M}$ commutes with arbitrary colimits in each factor.

\begin{lemma}\label{lem:inj_M}
Let $\msc{M}$ be a module category over a tensor category $\msc{A}$.  Then $\msc{M}$ has enough injectives and, for an arbitrary object $V$ in $\msc{A}$ and injective $I$ in $\msc{M}$, the product $I\ot V$ is injective in $\msc{M}$.
\end{lemma}

\begin{proof}
Similar to the proof of Lemma \ref{lem:inj_A}.
\end{proof}

In practice, the action map $\ot:\msc{M}\times\msc{A}\to \msc{M}$ will be exact in both factors.

\begin{lemma}\label{lem:biexact}
For any module category $\msc{M}$ over a tensor category $\msc{A}$, the action map $\ot:\msc{M}\times \msc{A}\to \msc{A}$ is exact in the $\msc{A}$ variable.  Furthermore, the action map $\ot$ is exact in the $\msc{M}$ variable if and only if the product $M\ot V$ of any nonzero objects $M$ in $\msc{M}$ and $V$ in $\msc{A}$ remains nonzero in $\msc{M}$.
\end{lemma}

We omit the proof, since the result is not used in this text.



Call a module category $\msc{M}$ \emph{pointed} if $\msc{M}$ comes equipped with a distinguished compact simple object $\1_\msc{M}$.  Equivalently, a pointed module category is a module category $\msc{M}$ which comes equipped with a choice of a right $\msc{A}$-module map
\[
u_{\msc{M}}:\msc{A}\to \msc{M}
\]
which sends $\1_\msc{A}$ to a compact simple object in $\msc{M}$.

\begin{remark}
By restricting to the compacts we recover the standard notions of tensor categories and module categories, as considered in \cite{egno15}.  We only work in the cocomplete setting because it is more convenient.
\end{remark}

\begin{remark}
The faithfullness condition of Lemma \ref{lem:biexact} holds for free whenever $\msc{A}$ admits a compact injective object.  We do not know if it holds in general, cf.\ \cite[Definition 7.3.1]{egno15}.
\end{remark}

\subsection{Tannakian reconstruction}

Suppose that $\msc{A}$ is a tensor category which comes equipped with a tensor functor to vector spaces, $F:\msc{A}\to \opn{Vect}$.  From such a pair $(\msc{A},F)$ one constructs a uniquely associated Hopf algebra $A$ and an equivalence of tensor categories $\opn{Corep}(A)\overset{\sim}\to \msc{A}$ which fits into a diagram
\[
\xymatrix{
\opn{Corep}(A)\ar[dr]_{\opn{forget}}\ar[rr]^{\sim}&  & \msc{A}\ar[dl]^F\\
 & \opn{Vect} &
}
\]
\cite[Theorems 2.2.8, 2.4.2]{schauenburg92}.  Indeed, this algebra is given by the so-called coendomorphisms
\[
A=\opn{Coend}(F)=\Hom_{\opn{Cont}/k}(\End_k(F),k).
\]
We refer the reader to \cite{schauenburg92} for a detailed presentation of this topic.

\subsection{Coideal subalgebras and module coalgebras}

Let $A$ be a Hopf algebra.  A right coideal subalgebra $B$ in $A$ is a subalgebra for which $\Delta(B)\subseteq B\ot A$.  A quotient module coalgebra $A\to C$ is a coalgebra quotient whose kernel is a right ideal in $A$.  In this way $C$ inherits a right $A$-action for which $\Delta(c\cdot a)=c_1\cdot a_1\ot c_2\cdot a_2$.  These structures, coideal subalgebras and quotient module coalgebras, are dual to each other when $A$ is finite-dimensional.
\par

At a categorical level, representations for a coideal subalgebra $\Rep(B)$ form a right $\Rep(A)$-module category, and corepresentations for a quotient module coalgebra $\opn{Corep}(C)$ form a right $\opn{Corep}(A)$-module category.  Furthermore, these module categories are pointed via the trivial representation for $B$, and the trivial corepresentation for $C$, respectively.

\subsection{Normality for pointed module categories}

Let $\msc{M}=(\msc{M},\1)$ be a pointed module category over a tensor category $\msc{A}$.  For the sake of specificity, let us say that $\msc{A}$ acts on the right of $\msc{M}$.
\par

Call an object $M$ in $\msc{M}$ trivial if $M$ admits a surjection $\1^{\oplus I}\to M$ from some additive power of the unit.

\begin{lemma}
\begin{enumerate}
\item Any trivial object is isomorphic to a (possibly infinite) additive power of $\1$.
\item Any quotient of a trivial object in $\msc{M}$ is trivial.
\item Any subobject of a trivial object is trivial.
\item Any object $M$ in $\msc{M}$ has a unique maximal trivial subobject $\opn{triv}(M)\subseteq M$.
\end{enumerate}
\end{lemma}

\begin{proof}
(1)--(3) All follow from the fact that any quotient of a semisimple object in an abelian category is semisimple, as is any subobject.  See for example \cite[Remark 2.2, Theorem 2.4]{lam91}.  (4) Take $\opn{triv}(M)=\opn{im}\big(\oplus_{f\in \Hom_{\msc{M}}(\1,M)}\1\to M)$.
\end{proof}

We consider the structure map $u_{\msc{M}}:\msc{A}\to \msc{M}$, and define the kernel of $u_{\msc{M}}$ to be
\[
\opn{ker}(u_{\msc{M}})=\left\{
\begin{array}{c}
\text{The full subcategory of $W$ in $\msc{A}$}\\
\text{for which $u_{\msc{M}}(W)$ is trivial in $\msc{M}$}
\end{array}\right\}.
\]
At this point we have standard notions of normality for tensor functors \cite[Definition 3.4]{bruguieresnatale11} (cf.\ \cite{etingofgelaki17}), which we now extend to the setting of pointed module categories.

\begin{definition}
For a pointed module category $\msc{M}$ over $\msc{A}$, we say the structure map $u_{\msc{M}}:\msc{A}\to \msc{M}$ is normal if the following properties hold:
\begin{enumerate}
\item[(a)] $u_{\msc{M}}$ is exact.
\item[(b)] The kernel of $u_{\msc{M}}$ is a tensor subcategory in $\msc{A}$.
\item[(c)] For any object $V$ in $\msc{A}$ there exists a subobject $\opn{triv}_{\msc{M}}(V)\subseteq V$ for which $u(\opn{triv}_{\msc{M}}(V))=\opn{triv}(u_{\msc{M}}(V))$.
\end{enumerate}
\end{definition}

We recall that exactness of the functor $u_{\msc{M}}=\1_{\msc{M}}\ot-$ holds for free in most practical situations.  See Lemma \ref{lem:biexact}.  Such exactness implies uniqueness of the subobject $\opn{triv}_{\msc{M}}(V)$ in $V$, and functoriality of the construction $V\mapsto \opn{triv}_{\msc{M}}(V)$.

\section{Quantum group basics}

Recall that $k$ is an algebraically closed field of characteristic $0$.

\subsection{Data for semisimple algebraic groups}

Fix a choice of a semisimple algebraic group $G$ over $k$, along with a choice of maximal torus $T\subseteq G$.  For such a pairing of $G$ and $T$ we have the associated root lattice $Q$, weight lattice $P$, and character lattice $X$.  These lattices fit into a sequence of inclusions $Q\subseteq X\subseteq P$.  We let $\Phi$ denote the collection of roots for $G$ and fix a choice of simple roots $\Delta$ in $\Phi$.  This choice of $\Delta$ determines a collection $\Phi^+$ of positive roots in $\Phi$.
\par

We let $(-,-):X\times X\to \mathbb{Q}$ denote the unique scaling of the Killing form so that $(\gamma,\gamma)=2$ at any short root $\gamma$.  We have
\[
(\gamma,\gamma)=d_{\gamma} \cdot 2
\]
at a general root $\gamma$, where $d_{\gamma}\in \{1,2,3\}$.  Let $\msc{W}=\langle\sigma_{\gamma}:\gamma\in \Phi\rangle$ denote the Weyl group for $G$.
\par

Throughout this work $G$ is always a semisimple algebraic group with specified data as above, and $\mfk{g}=\opn{Lie}(G)$.  We call $G$ \emph{almost-simple} if its Dynkin diagram is connected, and the lacing number for an almost-simple algebraic group is defined as
\[
\text{lacing number for G}:=\ \opn{max}\left\{-\frac{2(\alpha,\beta)}{(\alpha,\alpha)}:\ \alpha,\beta\ \in\Delta,\  \alpha\neq \beta\right\}=\max\{d_\gamma:\gamma\in \Phi\}.
\]
This number is $1$, $2$, and $3$, in types $\{A,D,E\}$, $\{B,C,F\}$, and $G_2$ respectively.

\subsection{Data for quantum groups}
\label{sect:q_data}

Consider the collection of connected components $\pi_0(\Delta)$ in the Dynkin diagram for $G$, and for each $i\in \pi_0(\Delta)$ let $\Delta_i$ denote the associated component in $\Delta$.  Each component $\Delta_i$ determines an almost-simple subgroup $G_i$ in $G$, and the decomposition $\Delta=\Delta_1\amalg\dots\amalg\Delta_t$ corresponds to the decomposition of $G$ into almost-simple factors
\[
G=G_1\times\dots\times G_t.
\]

The category of quantum group representations is specified by the data of a semisimple algebraic group $G$, with fixed maximal torus and simple roots $\Delta$, and a choice of symmetric bicharacter on the weight lattice
\[
q:P\times P\to k^\times
\]
which exponentiates the Killing form.  Specifically, we require that the form $q$ is invariant under the action of the Weyl group and that $q(\lambda,\mu)=1$ whenever $(\lambda,\mu)=0$.  We call $q$ \emph{torsion} if $q$ takes values in the roots of unity $\opn{tors}(k^\times)$.
\par

Such a form $q$ determines unique scalars $q_i$ over each connected component $i\in \pi_0(\Delta)$ under which
\[
q(\alpha,\lambda)=q_i^{(\alpha,\lambda)}\ \ \text{whenever}\ \alpha\in \Delta_i,\ \lambda\in P.
\]
Indeed, if we let $\omega_\alpha\in P$ denote the fundamental weight associated to a given simple root $\alpha\in \Delta$, and define
\[
q_{\alpha}=q(\alpha,\omega_\alpha)\ \ \text{for simple}\ \alpha,
\]
we may take $q_i=q_\alpha$ at any short root $\alpha$ in $\Delta_i$.
\vspace{2mm}

\emph{We refer to such a bilinear form $q$ as a \emph{quantum parameter} for $G$.  We refer to the units $q_i$ over the various components $i\in \pi_0(\Delta)$ as the \emph{scalar parameters} for $G$ at $q$.}  When we speak of a quantum group at a $n$-th ``roots of unity", we mean a pairing of $G$ and $q$ for which the associated scalar parameters $q_i$ are all primitive $n$-th roots of $1$.
\vspace{2mm}

The values $q_\alpha$ determine a unique function $q_\ast:\Phi\to k^\times$ which is invariant under the action of the Weyl group and which takes the prescribed values $q_\alpha$ on simple roots.  Explicitly, $q_\gamma=q_i^{d_\gamma}$ for each $\gamma$ in the subsystem $\Phi_i=\Phi\cap \mathbb{Z}\Delta_i$.

Supposing the form $q$ is torsion, we define
\begin{equation}\label{eq:l_gamma}
l_\gamma=\opn{ord}\big(q(\gamma,\gamma)\big)=\opn{ord}(q^2(\gamma,-))\ \ \text{for each }\gamma\in \Phi,
\end{equation}
and for a given component $i\in \pi_0(\Delta)$ we take
\[
l_i=l_{\alpha}\ \ \text{where $\alpha$ is any short root in }\Delta_i.
\]

\vspace{2mm}
\emph{Throughout this text, all quantum parameters $q$ are assumed to be torsion!}
\vspace{2mm}

\subsection{The quantum enveloping algebra}
\label{sect:qea}

Suppose first that $G$ is almost-simple, i.e.\ that the Dynkin diagram for $G$ is connected.  We consider Lusztig's integral divided power algebra $U_v=U_v(\mfk{g})$ over $\mathbb{Z}[v,v^{-1}]$ \cite[\S 1.3, Theorem 6.7]{lusztig90II}, and define $U_q$ as the base change of $U_v$ along the algebra map $\mathbb{Z}[v,v^{-1}]\to k$, $v\mapsto q_i$, where $i$ labels the unique component in the Dynkin diagram $\Delta_i=\Delta$.  We might also write $U_q=U_{q_i}$ in this case.

When $G$ has multiple almost-simple factors $G=G_1\times\dots\times G_t$ we take
\[
U_q=U_{q_1}(\mfk{g}_1)\ot_k\dots\ot_k U_{q_t}(\mfk{g}_t).
\]
As one expects, $\mfk{g}_i=\opn{Lie}(G_i)$ in the above expression and the parameters $q_i$ are as in Section \ref{sect:q_data}.

We recall that $U_q$ has generators
\[
E^{(n)}_\alpha,\ F^{(n)}_\alpha,\ K_\alpha,\ \text{for simple $\alpha\in \Delta$ and }n\in\mathbb{Z}_{\geq 0},
\]
and we have the distinguished toral elements
\[
\binom{K_\alpha;0}{l_\alpha}=\prod_{i=1}^{l_\alpha}\left.\frac{K_\alpha v^{-i+1}-K_{\alpha}^{-1}v^{i-1}}{v^i-v^{-i}}\right|_{v=q_\alpha}
\]
which also appear with some frequency.  These generators satisfy the usual quantum Serre relations and we have a Hopf structure on $U_q$ which satisfies
\[
\Delta(E_\alpha)=E_\alpha\ot 1+K_\alpha\ot E_\alpha,\ \ \Delta(F_\alpha)=F_\alpha\ot K_\alpha^{-1}+1\ot F_\alpha,\ \ \Delta(K_\alpha)=K_\alpha\ot K_\alpha
\]
\cite[\S 1.1]{lusztig90II}.  We also consider the positive and negative subalgebras $U^+_q$ and $U^-_q$ in $U_q$ which are generated by the divided powers of the simple root vectors $E^{(n)}_\alpha$ and $F^{(n)}_\alpha$ respectively.
\par

\subsection{Vector space bases}

We enumerate the simple roots $\Delta=\{\alpha_1,\dots,\alpha_{\opn{rk}(G)}\}$ for $G$ and, for a reduced expression of the longest element in the Weyl group $w_0=s_r\dots s_1$, we obtain an enumeration of the positive roots as
\[
\Phi^+=\{\gamma_1,\dots,\gamma_r\},\ \ \gamma_i=s_t\dots s_{i+1}(\alpha_i).
\]
\par

As shown in \cite{lusztig90,lusztig90II,lusztig93}, the simple reflections $s_\alpha\in \msc{W}$ lift to algebra automorphisms $T_\alpha:U_q\to U_q$ which satisfy the relations of the type $\mfk{g}$ braid group $B_\mfk{g}$ \cite[\S 37.1, Ch 39]{lusztig93}.  Via applications of these braid group operators we obtain associated root vectors
\[
E_{\gamma_i}=T_t\dots T_{i+1}(E_i)\ \in\ U_q^+
\]
for each $\gamma\in \Phi^+$.  We obtain corresponding negative root vectors $F_{\gamma}=\phi(E_{\gamma})$ in $U^-_q$ via an application of the anti-algebra isomorphism $\phi:U_{q^{-1}}\to U_q$ which exchanges $E_\alpha$ and $F_\alpha$, and inverts $K_\alpha$.  We then have Lusztig's bases for $U_q^{\pm}$ via ordered monomials in the roots vectors
\[
\{E_{\gamma_1}^{(m_1)}\dots E_{\gamma_r}^{(m_r)}:m:\Phi^+\to \mathbb{Z}_{\geq 0}\}\ \ \text{and}\ \ \{F_{\gamma_1}^{(m_1)}\dots F_{\gamma_r}^{(m_r)}:m:\Phi^+\to \mathbb{Z}_{\geq 0}\}
\]
\cite[Proposition 41.1.3]{lusztig93} \cite[Proposition 6.7]{lusztig90II}.

\begin{lemma}
Consider a semisimple algebraic group $G$ at a torsion parameter $q$, and a positive root $\gamma\in \Phi^+$.  We have $l_\gamma=1$ if and only if $E_\gamma$ and $F_\gamma$ are non-nilpotent, and otherwise
\[
l_\gamma=\opn{nil.\!deg}(E_\gamma)=\opn{nil.\!deg}(F_\gamma).
\]
\end{lemma}

Here by $\opn{nil.\!deg}(x)$ we mean the smallest postive integer $m$ at which $x^m=0$.

\begin{proof}
Follows from the expressions
\[
E_\gamma^m=[m]_{q_\gamma}!\cdot E_\gamma^{(m)},\ \ F_\gamma^m=[m]_{q_\gamma}!\cdot F_\gamma^{(m)},
\]
and the fact that the $q_\gamma$-factorial $[m]_{q_\gamma}!$ \cite[\S 1.3.3]{lusztig93} vanishes if and only if $m\geq l_\gamma>1$.
\end{proof}



\begin{remark}
In this text we employ the specific braid group operators $T_\alpha:=T''_{\alpha,1}$ from \cite[\S 37.1]{lusztig93}.
\end{remark}

\subsection{Quantum group representations}

We let $\Rep(G_q)$ denote the tensor category of character graded representations for the quantum enveloping algebra $U_q$, in the precise sense of \cite[\S\ 2.4]{negron21} (cf.\ \cite{andersenpolowen91,andersenpolowen92}).  So, a given $G_q$-representation $V$ is graded by the character lattice $X$ for $G$, the divided powers $E^{(m)}_\alpha$ and $F^{(m)}_\alpha$ act on $V$ by homogeneous endomorphisms of degrees $m\cdot \alpha$ and $-m\cdot \alpha$ respectively, and the toral elements in $U_q$ act as the semisimple operators
\[
K_\alpha\cdot v= q(\alpha,\lambda)\!\ v\ \ \text{and}\ \ \binom{K_\alpha;0}{l_\alpha}\cdot v=\binom{\langle \alpha,\lambda\rangle}{l_\alpha}_{q_\alpha} v,\ \ \text{for }v\in V_\lambda.
\]
\par

The category $\Rep(G_q)$ is, equivalently, the category of integrable, unital representations for Lusztig's modified algebra $\Rep(G_q)=\Rep(\dot{U}_q)$ \cite[\S\ 31.2]{lusztig93}.

\begin{remark}
Up to isomorphism, the algebra $U_q$ only depends on the restriction of the form $q$ to $Q\times Q$, while the definition of the category $\Rep(G_q)$ depends on the restriction of $q$ to $Q\times X$.
\par

The restricted form $q|_{X\times X}$ provides solutions to certain equations which one needs to define the $R$-matrix (braiding) for $\Rep(G_q)$.  The existence of an extension of this form to the entire weight lattice allows us to extract the scalar parameters $q_i$ from $q$, and also provides a simply-connected form $G^{sc}_q$ for any $G_q$.
\end{remark}

\begin{remark}
Our use of bilinear forms, rather than roots of unity, loosely follows the presentation of \cite{gaitsgory21}.
\end{remark}

\subsection{Quantum function algebras}

The quantum function algebra is, by definition, the Hopf algebra one reconstructs from the category $\Rep(G_q)$ and its forgetful functor to vector spaces,
\[
\O(G_q):=\opn{Coend}(forget:\Rep(G_q)\to \opn{Vect}).
\]
Again we refer the reader to \cite{schauenburg92} for a thorough discussion of Tannakian reconstruction.

\subsection{Dominant weights and simple representations}

The following result is proved in \cite[Proposition 6.4]{lusztig89} at odd order $q$, and one can use the formula \cite[Lemma 32.1.2]{lusztig93} to extend the arguments of \cite{lusztig89} to general $q$.  See alternatively \cite[Propositions 31.2.7, 31.3.2]{lusztig93}.

\begin{proposition}[{\cite[Proposition 6.4]{lusztig89}}]\label{prop:dominant}
For each dominant weight $\lambda\in X^+$ there is a unique simple object $L(\lambda)$ in $\Rep(G_q)$ which is of highest weight $\lambda$.  Furthermore, the assignment
\[
X^+\to \opn{Irrep}(G_q),\ \ \lambda\mapsto L(\lambda),
\]
is a bijection.
\end{proposition}

\subsection{Collected notations}

Fix a semisimple algebraic group $G$ with a choice of maximal torus and simple roots.  We gather all of the notations from this section.

\begin{itemize}
\item $\Phi$ denotes the roots for $G$, $\Delta$ denotes the chosen base in $\Phi$, and $\Phi^+$ denotes the corresponding subset of positive roots in $\Phi$.  We identify $\Delta$ with the Dynkin diagram for $G$, via an abuse of notation.
\item $Q$, $P$, and $X$ denote the root, weight, and character lattices for $G$.
\item $(-,-)$ denotes the normalized Killing form on $P$, so that $(\alpha,\alpha)=2$ at any short root $\alpha$.  At a general root $\gamma\in \Phi$ we let $d_\gamma$ denote the relative length $d_\gamma=(\gamma,\gamma)/2$.
\item For any $i\in \pi_0(\Delta)$, $\Delta_i\subset \Delta$ denotes the corresponding connected component in the Dynkin diagram.  We let $G_i\subseteq G$ denote the associated almost-simple factor in $G$.
\item $q$ denotes a quantum parameter for $G$.  In particular, $q$ is a symmetric form on the weight lattice which exponentiates the Killing form.  We always assume $q$ is torsion, i.e.\ takes values in the roots of unity $\opn{tors}(k^\times)$.
\item For any simple root $\alpha$, with corresponding fundamental weight $\omega_\alpha$, we take $q_\alpha=q(\alpha,\omega_\alpha)$.  For a general positive root $\gamma$ we take $q_\gamma=q_\alpha$, where $\alpha$ is any simple root which is in the Weyl group orbit of $\gamma$.
\item For any $i\in \pi_0(\Delta)$, we take $q_i=q_\alpha$ where $\alpha$ is any short root in $\Delta_i$.  We note that $q_i$ does not independent of the choice of $\alpha$.  These $q_i$ are the \emph{scalar parameters} associated to $G$ at $q$.
\item For any $\gamma\in \Phi^+$ we take
\[
l_\gamma=\opn{ord}(q(\gamma,\gamma))=\opn{ord}(q^2(\gamma,-))=\left\{\begin{array}{cl}
\opn{nil.\!deg}(E_\gamma) & \text{if $E_\gamma$ is nilpotent}\\
1 & \text{otherwise.}\end{array}\right.
\]
\end{itemize}
Let us also list the quantum algebras which we'll encounter in the text, though one might ignore the interruption upon first reading.
\begin{itemize}
\item $U_q$ denotes Lusztig's divided power quantum enveloping algebra \cite{lusztig90,lusztig90II}.
\item $\dot{U}_q$ denotes Lusztig's modified enveloping algebra from \cite[Ch 31]{lusztig93}.
\item $\mfk{v}_q$ denotes the distinguished finite-dimensional Hopf subalgebra in $U_q$ from \cite[\S 8.2]{lusztig90II}, and $\dot{\mfk{u}}_q$ denotes the modified analog of $\mfk{v}_q$, as in \cite[\S 36.2.1]{lusztig93}.  See also Sections \ref{sect:vq1}, \ref{sect:vq2}, \ref{sect:dotuq}.
\item $\bar{u}_q$ denotes the ``smallest quantum algebra" from Sections \ref{sect:smallest}, \ref{sect:smallest2}.  There is also an intermediate Hopf algebra $v_q$ which appears in these sections, and which serves as a link between $\mfk{v}_q$ and $\bar{u}_q$.
\end{itemize}
As defined above, $\Rep(G_q)$ is the category of character graded representations for the divided power algebra $U_q$.

\section{Quantum Frobenius}

We recall Lusztig's construction of the quantum Frobenius functor from \cite{lusztig93}.  This functor specifies a dual group $G^\ast$ to $G$, a dual parameter $\varepsilon$ to $q$, and an exact monoidal embedding $\Rep(G^\ast_\varepsilon)\to \Rep(G_q)$.  The associated scalar parameters $\varepsilon_i$ for $\varepsilon$ take values $\pm 1$, so that $G^\ast_\varepsilon$ has representation theory which is approximately equivalent to the classical representation theory for $G^\ast$.
\par

The quantum Frobenius functor plays a fundamental role in the analyses of quantum group representations that follow.

\subsection{Quantum groups at quasi-classical parameters}

\begin{proposition}[{\cite[Proposition 33.2.3]{lusztig93}}]\label{prop:731}
Let $H$ be any semisimple algebraic group and $\varepsilon$ be a quantum parameter with all $\varepsilon_i=\pm 1$.  There is a linear equivalence of categories $\Rep(H_\varepsilon)\overset{\sim}\to \Rep(H)$.  In particular, $\Rep(H_\varepsilon)$ is semisimple.
\end{proposition}

Note that we have not claimed that $\Rep(H_\varepsilon)$ and $\Rep(H)$ are equivalent as tensor categories, at general $\varepsilon$.

\subsection{``Dual" lattices}

We take
\[
X^\ast:=\{\lambda\in X:q(\alpha,\lambda)^2=1\ \text{for all simple $\alpha$ with }l_\alpha>1\}
\]
\[
\hspace{2cm}=\big\{\lambda\in X:(\alpha,\lambda)\in l_i\mathbb{Z}\ \text{for each }\alpha\in \Delta_i\text{ and }i\in \pi_0(\Delta)\big\}.
\]
For the equality of these two expressions, note that $q(\alpha,\lambda)^2=q_i^{2(\alpha,\lambda)}$ whenever $\alpha\in \Delta_i$, so that $q(\alpha,\lambda)^2=1$ if and only if $(\alpha,\lambda)\in l_i\mathbb{Z}$.  For simple roots with $l_\alpha=1$ such vanishing of $q(\alpha,\lambda)^2$ holds for free, at all $\lambda$, by the definition of $l_\alpha$.
\par

We also consider the ``$l$-dualized" root lattice
\[
lQ:=\mathbb{Z}\cdot \{l_\alpha \alpha:\alpha\ \text{simple}\}.
\]
The following is essentially stated in \cite{lusztig93}.

\begin{lemma}\label{lem:Q_v_X}
$l Q$ is contained in $X^\ast$, and both of the sublattices $l Q$ and $X^\ast$ are stable under the action of the Weyl group on $P$.
\end{lemma}

\begin{proof}
For each pair of simple roots $\alpha$ and $\beta$ we have
\[
q(\alpha,l_\beta\beta)^2=q(\alpha,\beta)^{2l_\beta}=1,
\]
by the definition of $l_\beta$.  This establishes the containment $l Q\subseteq X^\ast$.
\par

As for stability under the action of the Weyl group, for the simple reflection $\sigma_\alpha$ associated to $\alpha\in \Delta_i$, and $\lambda\in X^\ast$, we have
\begin{equation}\label{eq:507}
\sigma_\alpha(\lambda)=\lambda-\frac{(\alpha,\lambda)}{d_\alpha}\alpha.
\end{equation}
Each pairing $(\alpha,\lambda)$ lies in $d_\alpha \mathbb{Z}$ and, by the definition of $X^\ast$, in $l_i\mathbb{Z}$ as well.
\par

When $d_\alpha\mid l_i$ we have $l_\alpha=l_i/d_\alpha$ so that $(\alpha,\lambda)/d_\alpha$ is seen to lie in $l_\alpha\mathbb{Z}$.  When $d_\alpha\nmid l_i$ we have $l_\alpha=l_i$ and $(\alpha,\lambda)\in \opn{lcm}(d_\alpha,l_i)\mathbb{Z}=d_\alpha l_i\mathbb{Z}=d_\alpha l_\alpha \mathbb{Z}$, giving again $(\alpha,\lambda)/d_\alpha\in l_\alpha \mathbb{Z}$.  So from the expression \eqref{eq:507} we find
\begin{equation}\label{eq:516}
\sigma_\alpha(\lambda)\ \in\ \lambda+l_\alpha\mathbb{Z}\cdot \alpha\ \subseteq\ \lambda+l Q.
\end{equation}
This expression lies in $lQ$ when $\lambda$ lies in $lQ$, and lies in $X^\ast$ when $\lambda$ lies in $X^\ast$.  So we observe stability under the action of the Weyl group.
\end{proof}

\subsection{Quasi-classical representations in $\Rep(G_q)$}

\begin{proposition}[{\cite[Proposition 35.3.2]{lusztig93}}]\label{prop:674}
A simple representation $L(\lambda)$ in $\Rep(G_q)$ is annihilated by all simple root vectors $E_\alpha$ and $F_\alpha$ with $l_\alpha>1$ if and only if $\lambda\in X^\ast$.  Furthermore, in this case $L(\lambda)$ is graded by the sublattice $X^\ast$ in $X$.
\end{proposition}

\begin{proof}
If $\lambda\in X^\ast$ then $L(\lambda)$ has grading in $X^\ast$, and is annihilated by all such $E_\alpha$ and $F_\alpha$, directly by \cite[Proposition 35.3.1 (a) and (b)]{lusztig93}.  If $\lambda\notin X^\ast$ consider a simple root $\alpha$ with $l_\alpha>1$ and $q(\alpha,\lambda)^2\neq 1$, and $v\in L(\lambda)$ of highest weight to find
\[
[E_\alpha, F_\alpha]\cdot v =\frac{q(\alpha,\lambda)-q(\alpha,\lambda)^{-1}}{q_\alpha-q_\alpha^{-1}}\cdot v=q(\alpha,\lambda)^{-1}\left(\frac{q(\alpha,\lambda)^2-1}{q_\alpha-q_\alpha^{-1}}\right)\cdot v\neq 0.
\]
This equation implies $F_\alpha\cdot v\neq 0$ and $E_\alpha (F_\alpha\cdot v)\neq 0$.  Hence $L(\lambda)$ is not annihilated by the prescribed root vectors.
\end{proof}

\begin{lemma}
For a $G_q$-representation $V$, the following are equivalent:
\begin{enumerate}
\item[(a)] $V$ is annihilated by all root vectors $E_\alpha$ and $F_\alpha$, with $\alpha$ simple and $l_\alpha>1$.
\item[(b)] $V$ is annihilated by all root vectors $E_\gamma$ and $F_\gamma$, with $\gamma\in \Phi^+$ and $l_\gamma>1$.
\end{enumerate}
\end{lemma}

\begin{proof}
The implication (b) $\Rightarrow$ (a) is trivial.  So we consider (a) $\Rightarrow$ (b).  Let $I$ be the ideal in $U_q$ generated by the simple root vectors $E_\alpha$ and $F_\alpha$ with $l_\alpha>1$.  It suffices to show that $I$ is stable under the actions of the simple braid group operators $T_\beta$.  Indeed for any positive root $\gamma$ with $l_\gamma>1$ we can find a corresponding simple root $\alpha$, with $l_\alpha>1$, and a braid group operator $T$ for which $T(E_\alpha)=E_\gamma$ and $T(F_\alpha)=F_\gamma$.  However, such stability is clear from the explicit expressions of the $T_\beta$ given in \cite[\S 37.1.3]{lusztig93}.
\end{proof}

We call a representation $V$ in $\Rep(G_q)$ \emph{quasi-classical} if $V$ is annihilated by all simple $E_\alpha$ and $F_\alpha$ with $l_\alpha>1$.  Note that this class of representations is closed under taking subobjects and quotients.  We let $\msc{E}_q\subseteq \Rep(G_q)$ denote the full subcategory of quasi-classical representations.

\begin{lemma}\label{lem:696}
The subcategory $\msc{E}_q\subseteq \Rep G_q$ is a tensor subcategory in $\Rep G_q$, and all objects in $\msc{E}_q$ are graded by the sublattice $X^\ast$.
\end{lemma}

\begin{proof}
The fact that $\msc{E}_q$ is a tensor subcategory in $\Rep G_q$ follows from the fact that all simple root vectors $E_\alpha$ and $F_\alpha$ are skew primitive.  As for the claim about the grading, if a representation $V$ is annihilated by all $E_\alpha$ and $F_\alpha$ with $l_\alpha>1$ then all of its composition factors are of the form $L(\lambda)$ with $\lambda\in X^\ast$, by Proposition \ref{prop:674}.  It follows that $V$ is graded by $X^\ast$, since all of its composition factors are graded by $X^\ast$.
\end{proof}

\begin{remark}
The category $\msc{E}_q$ is also closed under the formation of extensions in $\Rep(G_q)$, provided all $l_i>2$, or provided $G$ has no factors of types $A_1$ or $B_n$ at $l_i=2$.  Indeed, in this case $\msc{E}_q$ is precisely the subcategory of $X^\ast$-graded objects in $\Rep(G_q)$.
\end{remark}

\subsection{The dual group $G^\ast$}

In \cite{lusztig90,lusztig90II,lusztig93} Lusztig considers a dual group $G^\ast$ to $G$ which is specified by the following data: The character lattice for $G^\ast$ is $X^\ast$, and the simple roots are
\[
\Delta^\ast:=\{\alpha^\ast:\alpha\in \Delta\}\ \subseteq\ l Q
\]
where $\alpha^\ast:=l_\alpha \alpha$.  The Cartan integers are given by
\[
\langle \alpha^\ast,\beta^\ast\rangle_l := \frac{2 l_\beta}{l_\alpha} \frac{ (\alpha,\beta)}{(\alpha,\alpha)}
\]
\cite[\S\ 2.2.5]{lusztig93}.  The group $G^\ast$ decomposes into almost-simple factors which correspond bijectively to the almost-simple factors for $G$,
\[
G=G_1\times\dots\times G_t\ \Rightarrow\ G^\ast=G^\ast_1\times\dots \times G^\ast_t.
\]
\par

One sees from the definition of the Cartan integers that the Killing form on the root lattice $l Q$ for $G^\ast$ is some rescaling of the form $(-,-)|_{l Q\times l Q}$.  We have $(\alpha^\ast,\alpha^\ast)=l_\alpha^2(\alpha,\alpha)$ and, for
\[
l'_i=\opn{min}\{l_\beta:\beta\in \Delta_i\},
\]
we rescale by $1/l_il_i'$ on each simple factor $\mfk{g}_i^\ast$ to obtain the properly normalized form $(-,-)_l$.  At each simple root $\alpha\in \Delta_i$ we now have
\[
\frac{1}{2}(\alpha^\ast,\alpha^\ast)_l=\frac{l_\alpha^2}{l_il_i'}d_\alpha=\left\{\begin{array}{ll}
d_\alpha &\text{if }l_i/l_i'=1\\
2 & \text{if $\alpha$ short and }l_i/l'_i=2\\
3 & \text{if $\alpha$ short and }l_i/l'_i=3\\
1 & \text{if $\alpha$ long and }l_i/l'_i=2\text{ or }3
\end{array}\right. .
\]
We therefore observe the fundamental weights
\[
\omega_{\alpha^\ast}=l_\alpha\cdot \omega_\alpha
\]
and weight lattice $P^\ast=\mathbb{Z}\cdot\{l_\alpha\omega_\alpha:\alpha\ \text{simple}\}$ for $G^\ast$.  From the above expressions we find the following.

\begin{lemma}
When the lacing number for $G_i$ is coprime to $l_i$, the dual group $G^\ast_i$ is of the same Dynkin type as $G_i$.  When the lacing number for $G_i$ divides $l_i$, the dual group $G^\ast_i$ is of Langlands dual Dynkin type to $G_i$.
\end{lemma}

\subsection{Lusztig's quantum Frobenius}
\label{sect:frob}

Take $G^\ast$ and $P^\ast$ as above.  Consider the form $\varepsilon$ on $P^\ast$ defined simply by restricting $q$,
\[
\varepsilon(\lambda,\mu):=q(\lambda,\mu)\ \ \text{for}\ \ \lambda,\mu\in P^\ast.
\]
Since the Killing form on $P^\ast$ is just a rescaling of the Killing form on the ambient lattice $P$, the form $\varepsilon$ exponentiates the Killing form for $G^\ast$ and so provides a quantum parameter for the dual group.  We have
\[
\varepsilon(\alpha^\ast,\alpha^\ast)=q(\alpha,\omega_\alpha)^{2l_\alpha^2}=1
\]
at each simple root, which forces $\varepsilon_i=\pm 1$ on each almost-simple factor $G_i^\ast$.

\begin{remark}
Scalar parameters of non-trivial value $\varepsilon_i=-1$ can occur, for example, when the corresponding parameter $q_i$ for $G_i$ is a $2p$-th root of $1$ with $p$ odd.  This can also occurs in types $B_n$, $C_n$, and $F_4$ when $q_i$ is a $4p$-th root of $1$ with $p$ odd, and in type $G_2$ when $q_i$ is a $6p$-th root of $1$ with $p$ odd.
\end{remark}

We now let $\dot{U}^\ast_{\varepsilon}$ denote the modified quantum enveloping algebra for $G^\ast$ at $\varepsilon$ \cite[\S 31.2]{lusztig93}, and let
\[
e^{(m)}_\alpha 1_\mu,\ f^{(m)}_\alpha 1_\mu\ \in \dot{U}^\ast_\varepsilon
\]
denote the standard generators, for simple $\alpha$ and $\mu\in X^\ast$.  By \cite[Theorem 35.1.9, \S\ 35.5.2]{lusztig93} (see also \cite{lentner16}) there is a surjective map of non-unital algebras
\[
fr:\dot{U}_q\to \dot{U}^\ast_\varepsilon
\]
which maps the powers of the simple root vectors as
\[
E^{(n)}_\alpha 1_\lambda\mapsto \left\{\begin{array}{ll}
e^{(n/l_\alpha)}_\alpha 1_\lambda &\text{if $l_\alpha |n$ and }\lambda \in X^\ast\vspace{2mm}\\
0 &\text{ otherwise}\end{array}\right.
\]
and
\[
F^{(n)}_\alpha 1_\lambda\mapsto \left\{\begin{array}{ll}
f^{(n/l_\alpha)}_\alpha 1_\lambda &\text{if $l_\alpha |n$ and }\lambda \in X^\ast\vspace{2mm}\\
0 & \text{ otherwise}.\end{array}\right.
\]
\par

The map $fr$ is compatible with the Hopf structures on $\dot{U}_q$ and $\dot{U}^\ast_\varepsilon$, and restricting along $fr$ provides an embedding of tensor categories
\begin{equation}
\opn{Fr}:=\res_{fr}:\Rep(G^\ast_{\varepsilon})\to \Rep(G_q).
\end{equation}
The tensor compatibility $\opn{Fr}(V)\ot \opn{Fr}(W)\to \opn{Fr}(V\ot W)$ is given by the identity.

\subsection{Quasi-classical representations via quantum Frobenius}

\begin{theorem}[{\cite{lusztig93}}]\label{thm:LQF}
The quantum Frobenius functor $\opn{Fr}:\Rep(G^\ast_\varepsilon)\to \Rep(G_q)$ restricts to an equivalence onto the subcategory $\msc{E}_q$ of quasi-classical representations in $\Rep(G_q)$.
\end{theorem}

The proof is essentially the same as that of \cite[Corollary 35.3.3]{lusztig93}.  We repeat it for the sake of completeness.

\begin{proof}
We first note that the image of $\opn{Fr}$ is contained in $\msc{E}_q$, directly by the definitions.  Now, it is shown in \cite[\S\ 35.2.3]{lusztig93} that the generators $E^{(l_\alpha)}_\alpha$ in $U^+_q$ satisfy the $\varepsilon$-Serre relations of the $e_\alpha$ in $U^+_\varepsilon$.  The $F^{(l_\alpha)}_\alpha$ similarly satisfy the relations of the $f_\alpha$.  So any $V$ in $\msc{E}_q$ comes equipped with actions of $U_{\varepsilon}^+$ and $U_{\varepsilon}^-$ and is appropriately graded by $X^\ast$.  The commutativity relations \cite[6.5 (a2)]{lusztig90II}, in conjunction with the fact that $V$ is annihilated by all $E_\alpha$ with $l_\alpha>1$, implies that
\[
[E^{(l_\alpha)}_\alpha,F^{(l_\beta)}_\beta]\cdot v=\delta_{\alpha\beta}\binom{K_\alpha;0}{l_\alpha}\cdot v=\delta_{\alpha,\beta}\binom{\langle \alpha,\lambda\rangle}{l_\alpha}_{q_\alpha}\cdot v\\
\]
\[
=\delta_{\alpha\beta}\left(q_\alpha^{(\frac{\langle \alpha,\lambda\rangle}{l_\alpha}+1)l^2_\alpha}\frac{\langle \alpha,\lambda\rangle}{l_\alpha}\cdot v\right)=\delta_{\alpha\beta}\left(q_\alpha^{l_\alpha^2(\langle \alpha^\ast,\lambda\rangle_l+1)}\langle\alpha^\ast,\lambda\rangle_l\cdot v\right)
\]
for arbitrary $\alpha, \beta \in \Delta$ \cite[Lemma 34.1.2]{lusztig93}.  We have $q_\alpha^{l_\alpha^2}=\varepsilon_{\alpha^\ast}$ so that the final expression reduces to
\[
[E^{(l_\alpha)}_\alpha,F^{(l_\beta)}_\beta]\cdot v=\delta_{\alpha\beta}\left(\varepsilon_{\alpha^\ast}^{(\langle \alpha^\ast,\lambda\rangle_l+1)}\langle \alpha^\ast, \lambda\rangle_l\cdot v\right)=\delta_{\alpha,\beta}\binom{\langle \alpha^\ast, \lambda\rangle_l}{1}_{\varepsilon_{\alpha^\ast}}\cdot v.
\]
This information tells us that the $E^{(l_\alpha)}_\alpha$ and $F^{(l_\alpha)}_\alpha$ actions on $V$ satisfy the necessarily relations to provide a $\dot{U}^\ast_{\varepsilon}$-action on $V$.  This defines a functor
\[
-|_{U^\ast_{\varepsilon}}:\msc{E}_q\to \Rep(G^\ast_{\varepsilon}).
\]
\par

Now, for $V$ in $\msc{E}_q$, the identity map provides a natural isomorphism of $G_q$-representations $\opn{Fr}(V|_{\dot{U}^\ast_{\varepsilon}})\overset{\sim}\to V$, and for any $W$ in $\Rep(G^\ast_\varepsilon)$ the identity again provides a natural isomorphism of $G^\ast_\varepsilon$-representations $W\overset{\sim}\to \opn{Fr}(W)|_{\dot{U}^\ast_{\varepsilon}}$.  So $\opn{Fr}:\Rep(G^\ast_\varepsilon)\to \Rep(G_q)$ restricts to an equivalence onto $\msc{E}_q$.
\end{proof}

The following is an immediate corollary to Theorem \ref{thm:LQF} and Proposition \ref{prop:731}.

\begin{corollary}
The category $\msc{E}_q$ of quasi-classical representations in $\Rep(G_q)$ is semisimple.
\end{corollary}

\section{Small quantum algebras above the lacing number}
\label{sect:sqa}

We cover the construction of the smallest quantum algebra $\bar{u}_q=\bar{u}(G_q)$ at almost-every $q$.  We address all pairings of $G$ and $q$ save for types $B_n$, $C_n$, and $F_4$ at a $4$-th root of $1$, and type $G_2$ at a $3$-rd, $4$-th, or $6$-th root of $1$.  These low order cases are discussed independently in Section \ref{sect:sqa_at_lacing} below.  We also explain how the category of representations $\Rep(\bar{u}_q)$ inherits the natural structure of a pointed module category over $\Rep(G_q)$.
\par

We begin by recalling Lusztig's finite-dimensional Hopf subalgebra in $U_q$, which we denote by $\mfk{v}_q$.  We view $\mfk{v}_q$ as associated to the Lie algebra $\mfk{g}=\opn{Lie}(G)$, rather then the group $G$ itself.  We then construct an intermediate Hopf algebra $v_q$ which admits highly structured maps from both $\mfk{v}_q$ and our ``smallest quantum algebra" $\bar{u}_q$.  We use $v_q$, along with the aforementioned maps, to translate results for $\mfk{v}_q$ to $\bar{u}_q$;
\[
\xymatrixrowsep{3mm}
\xymatrix{
	&v_q\\
\mfk{v}_q\ar[ur]\ar@{<-->}[rr] & & \bar{u}_q.\ar[ul]
}
\]

\subsection{Lusztig's finite-dimensional subalgebra}
\label{sect:vq1}

Consider a semisimple algebraic group $G$, and suppose that each $l_i$ is greater than the lacing number for its associated almost-simple factors in $G$.  We consider Lusztig's divided power algebra $U_q$, which we think of as a ``copy" of the quantum enveloping algebra used in the construction of $\Rep(G_q)$.   To distinguish between these two algebras we let $\msf{K}_\alpha$ denote the standard grouplikes in our new copy of $U_q$.
\par

We begin with Lusztig's finite-dimensional subalgebra in $U_q$ \cite{lusztig90,lusztig90II}.  We adopt the unorthodox notation $\mfk{v}_q$ for this algebra, and recall the definition
\[
\mfk{v}_q=\mfk{v}_q(\mfk{g})=\left\{
\begin{array}{c}
\text{The subalgebra in $U_q$ generated}\\
\text{by the elements $E_\alpha$, $F_\alpha$, and $\msf{K}_\alpha$}\\
\text{for simple roots }\alpha\in \Delta
\end{array}\right\}\ \subseteq\ U_q.
\]
The algebra $\mfk{v}_q$ is a finite-dimensional Hopf subalgebra in $U_q$ which has basis
\[
\mathbf{B}=\{E_{\gamma_1}^{m_1}\dots E_{\gamma_r}^{m_r}\cdot \msf{K}_\nu\cdot F_{\gamma_1}^{m'_1}\dots F_{\gamma_r}^{m'_r}:\nu\in Q,\ \text{\rm $\gamma_i\in \Phi^+$, and }m_i, m'_i<l_{\gamma_i}\},
\]
under some ordering of the positive roots.
\par

From this basis one determines the dimension of $\mfk{v}_q$ as
\[
\dim(\mfk{v}_q)=(\prod_{\alpha\in \Delta}2l_\alpha)(\prod_{\gamma\in \Phi^+}l_\gamma)^2.
\]
The first factor accounts for the grouplikes in $\mfk{v}_q$ and the second factors accounts for the positive and negative subalgebras $\mfk{v}^{\pm}_q$ in $\mfk{v}_q$ \cite[Theorem 8.3]{lusztig90II}.

\subsection{Small quantum algebras}
\label{sect:smallest}

Recall our root lattice $lQ=\mathbb{Z}\cdot\{l_\alpha\alpha:\alpha\in \Delta\}$ and character lattice $X^\ast$ for the dual group $G^\ast$, from Section \ref{sect:frob}.  We first consider the character groups
\[
A:=(X/2\cdot l Q)^\vee\ \ \text{and}\ \ \bar{A}:=(X/X^\ast)^\vee,
\]
and note that the inclusion $2\cdot l Q\subseteq X^\ast$ dualizes to an inclusion $\bar{A}\subseteq A$.

For any $\nu\in Q$ we have the associated character $K_\nu:X\to k^\ast$, $K_\nu(\lambda)=q(\lambda,\nu)$, which vanishes on $2\cdot l Q$ and hence defines an element in $A$.  We also have a natural action of $A$ on $\mfk{v}_q$ by the Hopf algebra automorphisms
\[
\xi\cdot E_\gamma=\xi(\gamma)E_\gamma,\ \ \xi\cdot F_\gamma=\xi(-\gamma)F_\gamma,\ \ \xi\cdot \msf{K}_\alpha=\msf{K}_\alpha,
\]
and can consider the new Hopf algebra $\mfk{v}_q\rtimes A$.
\par

We note that the differences $\msf{K}_\nu\cdot K_\nu^{-1}$ are central grouplike elements in $\mfk{v}_q\rtimes A$, and take the quotient
\[
v_q:= \frac{\mfk{v}_q\rtimes A}{(\msf{K}_\nu\cdot K_\nu^{-1}:\nu\in Q)}.
\]
From the basis for $\mfk{v}_q$ we see that $\mfk{v}_q\rtimes A$ is a symmetric bimodule over the subgroup algebra generated by the differences $\msf{K}_\nu\cdot K_\nu^{-1}$, which is furthermore free on the left and right independently.  We then have the induced basis,
\[
\text{Basis for }v_q=\{E_{\gamma_1}^{m_1}\dots E_{\gamma_r}^{m_r}\cdot \xi\cdot F_{\gamma_1}^{m'_1}\dots F_{\gamma_r}^{m'_r}:\xi\in A,\ m_i,m'_i<l_{\gamma_i}\ \text{for all }i\},
\]
and the triangular decomposition for $\mfk{v}_q$ induces a triangular decomposition for $v_q$,
\[
\mfk{v}^+_q\ot kA\ot \mfk{v}^-_q\overset{\cong}\to v_q.
\]
\par

In our new algebra $v_q$, we scale by the grouplikes to produce normalized root vectors $\mathbf{E}_\gamma=K_\gamma E_\gamma$ for each $\gamma\in \Phi^+$.

\begin{definition}
Let $G$ be a semisimple algebraic group, and $q$ be a quantum parameter at which each $l_i$ is greater than the lacing number for its associated almost-simple factor in $G$.  We define $\bar{u}_q$ as follows:
\[
\bar{u}_q=\bar{u}(G_q):=\left\{\begin{array}{c}
\text{The subalgebra in $v_q$ generated by the}\\
\text{root vectors $\mathbf{E}_\alpha$ and $F_\alpha$, for $\alpha\in \Delta$,}\\
\text{and the subgroup of characters $\bar{A}$ in $A$}
\end{array}\right\}.
\]
\end{definition}

We refer to $\bar{u}_q$, somewhat informally, as the \emph{smallest quantum algebra} for $G$ at $q$.  We consider the subalgebras $\bar{u}^{+}_q$ and $\bar{u}^-_q$ in $\bar{u}_q$ generated by the positive and negative root vectors $\mathbf{E}_\alpha$ and $F_\alpha$ respectively.

\begin{lemma}\label{lem:973}
Under some ordering of the positive roots $\{\gamma_1,\dots,\gamma_r\}=\Phi^+$, the collections
\[
\{\mathbf{E}_{\gamma_1}^{m_1}\dots \mathbf{E}_{\gamma_r}^{m_r}:m_i<l_{\gamma_i}\ \text{for all }i\}\ \ \text{and}\ \ \{F_{\gamma_1}^{m'_1}\dots F_{\gamma_r}^{m'_r}:m'_i<l_{\gamma_i}\ \text{for all }i\}
\]
provide bases for $\bar{u}^+_q$ and $\bar{u}^-_q$ respectively.
\end{lemma}

\begin{proof}
Take $v^{\pm}_q=\bar{u}_q^{\pm}$, interpreted as subalgebras in the ambient algebra $v_q$.  For $v^+_q$ we note that both of the subalgebras $v^+_q,\mfk{v}^+_q\subseteq v_q$ are naturally graded by the root lattice $Q$, and we have linear isomorphisms
\[
K_\nu\cdot-:(\mfk{v}_q^+)_\nu\overset{\sim}\to (v^+_q)_\nu
\]
at all $\nu\in Q$.  So we deduce the claimed basis for $v^+_q$ from that of $\mfk{v}^+_q$.  For $\bar{u}_q^-$, we simply have $\bar{u}^-_q=v^-_q=\mfk{v}^-_q$.
\end{proof}

From the triangular decomposition for $v_q$ we deduce a triangular decompositions for $\bar{u}_q$,
\begin{equation}\label{eq:triang_baruq}
\bar{u}^+_q\ot k\bar{A}\ot \bar{u}^-_q\overset{\cong}\to \bar{u}_q.
\end{equation}
The above isomorphism is given explicitly by including the relevant algebras into $\bar{u}_q$, then composing with multiplication.  From the decomposition \eqref{eq:triang_baruq}, and the bases of  Lemma \ref{lem:973}, we calculate the dimension of $\bar{u}_q$.

\begin{corollary}\label{cor:dims}
\[
\dim(\bar{u}_q)=[X:X^\ast]\cdot (\prod_{\gamma\in \Phi^+}l_\gamma)^2.
\]
\end{corollary}


\begin{remark}
Our rescaling of the generators $E_\gamma\mapsto \mbf{E}_\gamma$ is adapted from Arkhipov-Gaitsgory \cite{arkhipovgaitsgory03}.
\end{remark}

\subsection{Tensor structures}
\label{sect:res_ten}

The Hopf structure on $v_q$ is given explicitly by the formulas $\Delta(\xi)=\xi\ot\xi$, for $\xi\in A$, and
\begin{equation}\label{eq:858}
\Delta(E_\alpha)=E_\alpha\ot 1+K_\alpha\ot E_\alpha,\ \ \Delta(F_\alpha)=F_\alpha\ot K_{-\alpha}+1\ot F_\alpha,
\end{equation}
for each $\alpha\in \Delta$.

\begin{lemma}
The algebra $\bar{u}_q$ is a right coideal subalgebra in $v_q$.
\end{lemma}

\begin{proof}
We have $\Delta(F_\alpha)\in \bar{u}_q\ot v_q$ and $\Delta(\xi)=\xi\ot \xi\in \bar{u}_q\ot v_q$ for all $\xi\in \bar{A}$, by direct inspection.  For any $\lambda\in X^\ast$ we have $q^2(\alpha,\lambda)=1$ at all $\alpha\in \Delta$, so that $K^2_\alpha:X\to k^\times$ vanishes on $X^\ast$.  Hence each character $K^2_\alpha$ lies in $\bar{A}\subseteq \bar{u}_q$.  So we have
\[
\Delta(\mathbf{E}_\alpha)=\mathbf{E}_\alpha\ot K_\alpha+K_\alpha^2\ot \mathbf{E}_\alpha\ \in\ \bar{u}_q\ot v_q,
\]
and see that $\bar{u}_q$ is a right coideal subalgebra in $v_q$.
\end{proof}

Via the actions of $\mathbf{E}_\alpha$, $F_\alpha$, and $\bar{A}$ on big quantum group representations, $\Rep(\bar{u}_q)$ inherits a natural right module category structure over $\Rep(G_q)$.

To elaborate, any $G_q$-representation carries a natural $X/X^\ast$-grading, via the $X$-grading and the projection $X\to X/X^\ast$.  Hence any $G_q$-representation carries a natural action of the character group $\bar{A}=(X/X^\ast)^\vee$.  Via this $\bar{A}$-action, the actions of $\mathbf{E}_\alpha=K_\alpha E_\alpha$ and $F_\alpha$ on $G_q$-representations, and the coproducts \eqref{eq:858} we obtain a natural action of $\bar{u}_q$ on products
\[
W\ot V\ \ \text{for $W$ in $\Rep(\bar{u}_q)$ and $V$ in }\Rep(G_q).
\]
In this way $\Rep(\bar{u}_q)$ is realized as a right module category over $\Rep(G_q)$.

Furthermore, the restriction functor
\[
\res^-_q:\Rep(G_q)\to \Rep(\bar{u}_q)
\]
is $\Rep(G_q)$-linear, with trivial compatibility $\res^-_q(V'\ot V)=\res^-_q(V')\ot V$, so that $\Rep(\bar{u}_q)$ becomes a \emph{pointed} module category over $\Rep(G_q)$.  The unit object in $\Rep(\bar{u}_q)$ is given by the trivial representation $\1_{\Rep(\bar{u}_q)}=k$.

We find in Section \ref{sect:normal} below that the restriction functor $\res^-_q:\Rep(G_q)\to \Rep(\bar{u}_q)$ is in fact normal, and has kernel equal to the subcategory of quasi-classical representations in $\Rep(G_q)$.

\section{Small quantum algebras below the lacing number}
\label{sect:sqa_at_lacing}

We extend our construction of $\bar{u}_q$ at not-too-small order $q$ to arbitrary $q$.  Here we replace Lusztig's usual small quantum algebra, from Section \ref{sect:vq1}, with a slight modification which already appears in \cite[\S 8.2]{lusztig93} and was studies extensively by Lentner \cite{lentner16}.
\par

Formally speaking, the materials of this section recover (and extend) all of our findings from Section \ref{sect:sqa}.  However, our presentation becomes much more delicate as we must account for degenerations of certain ``linkings" between root vectors in $U_q$ at low order $q$.

\begin{remark}
The reader who is happy to work away from extraordinarily small order parameters can safely skip this section, though they should heed the change in notation from $\Delta$ to $\Delta_l$ (see Section \ref{sect:delta_l}).
\end{remark}

\subsection{Finite-dimensional subalgebras at small order $q$}
\label{sect:vq2}

We consider the quantum enveloping algebra $U_q$, now at arbitrary $q$.  We again think of $U_q$ as a ``copy" of our original enveloping algebra from Section \ref{sect:qea}, and let $\msf{K}_\nu$ denote the grouplike element in $U_q$ associated to $\nu\in Q$.  In \cite[\S 8.1, 8.2]{lusztig90II} Lusztig introduces the following subalgebra in $U_q$,\footnote{It is not obvious that our algebra $\mfk{v}_q$ agrees with Lusztig's algebra from \cite[\S 8.1, 8.2]{lusztig90II}.  However, such an identification follows from Propositions \ref{prop:863} and \ref{prop:865} below in conjunction with stability under the action of the braid group.}
\begin{equation}\label{eq:mfkv_q}
\mfk{v}_q=\mfk{v}_q(\mfk{g})=\left\{\begin{array}{c}
\text{The smallest subalgebra in $U_q$ which contains all}\\
\text{of the simple generators $E_\alpha$ and $F_\alpha$ with $l_\alpha>1$,}\\
\text{contains all the grouplikes $\msf{K}_\nu$ for $\nu\in Q$, and }\\
\text{which is stable under the braid group action on $U_q$}
\end{array}\right\}.
\end{equation}
\par

We note that $\mfk{v}_q$ is \emph{not} generated by the simple root vectors $E_\alpha$ and $F_\alpha$ in general.  For example, in type $B_2$ at a $4$-th root of $1$ we have $l_\beta=1$ at the long simple root $\beta$, so that neither $E_\beta$ nor $F_\beta$ appears in $\mfk{v}_q$.  However, the root vectors $E_{\beta+\alpha}=E_\beta E_\alpha+E_\alpha E_\beta$ and $F_{\beta+\alpha}$ do appear in $\mfk{v}_q$ in this case.
\par

Lentner establishes the following structural result.

\begin{proposition}[{\cite[Theorems 5.2 \& 5.4]{lentner16}}]\label{prop:863}
Suppose that $G$ is almost-simple, and take $l=l_i$ for the unique component $\Delta_i=\Delta$.  At an arbitrary quantum parameter $q$ we have a triangular decomposition
\[
\mfk{v}_q^+\ot k[\msf{K}_\nu:\nu\in Q]\ot \mfk{v}_q^-\overset{\sim}\to \mfk{v}_q
\]
and corresponding linear basis
\begin{equation}\label{eq:910}
\mathbf{B}=\{E_{\gamma_1}^{m_1}\dots E_{\gamma_r}^{m_r}\cdot \msf{K}_{\nu}\cdot F_{\gamma_1}^{m'_1}\dots F_{\gamma_r}^{m'_r}:\nu\in Q\ \text{\rm and}\ m_i,\ m'_i<l_{\gamma_i}\}
\end{equation}
for $\mfk{v}_q$, where $\{\gamma_1,\dots,\gamma_r\}=\Phi^+$ is some ordering of the positive roots.
\par

In the specific cases where $l$ is equal to the lacing number for $G$, the algebra $\mfk{v}_q$ admits the following description:
\begin{itemize}
\item When $l=1$, $\mfk{v}_q=k[K_\nu:\nu\in Q]$.
\item In type $C_n$ at $l=2$, $\mfk{v}_q$ is generated by the simple root vectors $E_\alpha$ and $F_\alpha$ with $\alpha$ short, the grouplikes $K_\nu$, and
\begin{equation}\label{eq:926}
E_{\beta+\alpha'}:=(E_{\beta} E_{\alpha'}+E_{\alpha'}E_{\beta}),\ \ F_{\beta+\alpha'}:=-(F_{\beta} F_{\alpha'}+F_{\alpha'}F_{\beta}),
\end{equation}
where $\beta$ is the long simple root and $\alpha'$ is its unique short neighbor.
\item In type $F_4$ at $l=2$, list the roots as $\{\beta_4,\beta_3,\alpha_2,\alpha_1\}$ where the $\alpha_i$ are short and $(\beta_3,\alpha_2)=-2$.  Then $\mfk{v}_q$ is generated by the short vectors $E_{\alpha_i}$ and $F_{\alpha_i}$, the grouplikes $K_\nu$, $E_{\beta_3+\alpha_2}$ and $F_{\beta_3+\alpha_2}$ defined as in \eqref{eq:926}, and
\[
E_{\beta_4+\beta_3+\alpha_2}=(E_{\beta_4}E_{\beta_3+\alpha_2}+E_{\beta_3+\alpha_2}E_{\beta_4}),\ \ F_{\beta_4+\beta_3+\alpha_2}=-(F_{\beta_4}F_{\beta_3+\alpha_2}+F_{\beta_3+\alpha_2}F_{\beta_4}).
\]
\item In type $B_n$ at $l=2$, list the roots as $\{\beta_n,\dots,\beta_{2},\alpha_1\}$, where consecutive roots are neighbors and $\alpha_1$ is short.  Then $\mfk{v}_q$ is generated by $E_{\alpha_1}$ and $F_{\alpha_1}$, the grouplikes $K_\nu$, and the root vectors $E_{\beta_r+\cdots \beta_2+\alpha_1}$ and $F_{\beta_r+\cdots \beta_2+\alpha_1}$ for all $r\geq 2$, which we define recursively as
\[
E_{\beta_{r+1}+\cdots \alpha_1}=(E_{\beta_{r+1}}E_{\beta_r+\cdots \alpha_1}+E_{\beta_r+\cdots \alpha_1}E_{\beta_{r+1}}),
\]
\[
F_{\beta_{r+1}+\cdots \alpha_1}=-(F_{\beta_{r+1}}F_{\beta_r+\cdots \alpha_1}+F_{\beta_r+\cdots \alpha_1}F_{\beta_{r+1}}).
\]
\item In type $G_2$ at $l=3$, $\mfk{v}_q$ is generated by $E_\alpha$ and $F_\alpha$ for short $\alpha$, all the grouplikes $K_\nu$, and the vectors
\[
E_{\beta+\alpha}=(E_\beta E_{\alpha}-q_\alpha^{3} E_\alpha E_\beta),\ \ F_{\beta+\alpha}=-(F_\beta F_\alpha-q_{\alpha}^{3}F_\beta F_\alpha).
\]
\end{itemize}
Furthermore, each of the generating $E_\gamma$ and $F_\gamma$ above are skew primitive in $\mfk{v}_q$.
\end{proposition}

\begin{proof}
Each anomalous generator $E_\gamma=E_{\beta+\gamma'}$ (resp.\ $F_{\gamma}=F_{\beta+\gamma'}$) is obtained from the preceding generator $E_{\gamma'}$ (resp.\ $F_{\gamma'}$) via an application of the braid group operator $T_\beta$.  So it is clear that these elements live in $\mfk{v}_q$.  In \cite{lentner16} it is shown that these generators $E_\gamma$ and $F_\gamma$ are also skew primitive in $\mfk{v}_q$, and that they generate a finite-dimensional subalgebra in $U_q$ with the prescribed basis.  So, all that is left to show is that the subalgebra generated by the prescribed elements is all of $\mfk{v}_q$.  Equivalently, we must show that this subalgebra is stable under the braid group action on $U_q$.  We verify such stability via direct computation in Appendix \ref{sect:B_g}.
\end{proof}

\begin{remark}
Though the algebra $\mfk{v}_q$ already appears in \cite[\S 8.2]{lusztig90II}, Lusztig immediately relocates to the case of an odd order scalar parameter.  Such low order cases are also avoided in the text \cite{lusztig93}.  So we rely on, and emphasize, the more extensive analysis of Lentner \cite{lentner16} in this section.
\end{remark}

One should note that, in the above basis, we have $l_\gamma=1$ whenever $\gamma\in \Phi^+$ is long and $l$ is equal to the lacing number for $G$.  Hence many root vectors do not contribute to the basis $\mathbf{B}$ at low order $q$.  We also highlight the validity of the expression \eqref{eq:910} at \emph{all} quantum parameters (cf.\ Section \ref{sect:vq1}).
\par

The final case of $G_2$ at a $4$-th root of unity is also interesting.  Here not all root vectors lie in the subalgebra of $U_q$ generated by the simple root vectors, though all root vectors $E_\gamma$ and $F_\gamma$ do appear in $\mfk{v}_q$.

\begin{proposition}[{\cite[Theorem 5.4]{lentner16}}]\label{prop:865}
Consider $G$ of type $G_2$ at $l=2$, and let $\beta$ and $\alpha$ be the long and short simple roots respectively.  The algebra $\mfk{v}_q$ is generated by the grouplikes $K_\nu$ and the elements
\[
E_\alpha,\ E_\beta,\ E_{2\alpha+\beta},\ \ F_\alpha,\ F_\beta,\ F_{2\alpha+\beta},
\]
where
\[
E_{2\alpha+\beta}=E_\alpha^{(2)}E_\beta+E_\beta E_\alpha^{(2)}+E_\alpha E_\beta E_\alpha\ \ \text{and}\ \ 
F_{2\alpha+\beta}=F_\alpha^{(2)}F_\beta+F_\beta F_\alpha^{(2)}+F_\alpha F_\beta F_\alpha.
\]
Furthermore, these generating $E$'s and $F$'s are skew primitive.
\end{proposition}

\begin{proof}
We have
\[
E_{2\alpha+\beta}=T_\beta T_\alpha T_\beta(E_\alpha)\ \ \text{and}\ \ F_{2\alpha+\beta}=T_\beta T_\alpha T_\beta(F_\alpha)
\]
\cite[\S 39.2.2]{lusztig93}.  Lenter proves skew primitivity of these elements in $\mfk{v}_q$, and that the subalgebra in $U_q$ which is generated by these elements has the prescribed basis \cite[Theorem 5.4]{lentner16}.  So, again, we need only verify that the subalgebra generated by the specified elements $E_\gamma$, $F_\gamma$, and $K_\nu$ is stable under the braid group action.  We verify such stability in Appendix \ref{sect:B_g}.
\end{proof}

For general $G$, $\mfk{v}_q$ decomposes into a product
\[
\mfk{v}_q=\mfk{v}_{q}(\mfk{g}_1)\ot\dots\ot \mfk{v}_{q}(\mfk{g}_t)
\]
in accordance with the decomposition of $G$ into almost-simple factors.  So Propositions \ref{prop:863} and \ref{prop:865} provide an explicit description of $\mfk{v}_q$ in all cases.
\par

Though there are many more interesting things to say about this algebra, and in particular about its positive and negative subalgebras $\mfk{v}^{\pm}_q$ at low order $q$, we close our introduction to $\mfk{v}_q$ here and invite the reader to see \cite{lentner16} for more details.

\subsection{Distinguished root vectors at low order $q$}\label{sect:delta_l}

It will be convenient to define a subset of positive roots $\Delta_l\subseteq \Phi^+$ which labels the skew-primitive generators in $\mfk{v}_q$ at arbitrary $q$.  Though the definition of $\Delta_l$ is implicit in the statements of Propositions \ref{prop:863} and \ref{prop:865} above, let us take a moment to identify this subset clearly.
\par

First, we decompose $G$ into its almost-simple factors $G=G_1\times\dots\times G_t$, each of which has its associated simple roots $\Delta_i\subseteq \Delta$.  For each index $i\in \pi_0(\Delta)$ we consider a new set of ``simple roots" $\Delta_{l_i}$, and then define $\Delta_l=\coprod_i \Delta_{l_i}\ \subseteq\ \Phi^+$.

In the case where $l_i$ is greater than the lacing number for $G_i$ we just take $\Delta_{l_i}=\Delta$.  In the case where $l_i$ is equal to the lacing number for $G_i$, we first enumerate the simple roots as
\[
\Delta_i:=\{\beta_n,\dots,\beta_{r+1},\alpha_r,\dots,\alpha_1\}
\]
with the $\alpha_i$ short, the $\beta_j$ long, and all consecutive roots neighbors.  We take
\begin{itemize}
\item $\Delta_{l_i}=\{\beta_n+\dots +\beta_2+\alpha_1,\dots,\beta_2+\alpha_1,\alpha_1\}$ when $G_i$ is of type $B_n$ at $l_i=2$.\vspace{2mm}
\item $\Delta_{l_i}=\{\beta_n+\alpha_{n-1},\alpha_{n-1},\dots,\alpha_1\}$ in type $C_n$ at $l_i=2$.\vspace{2mm}
\item $\Delta_{l_i}=\{\beta_4+\beta_3+\alpha_2,\beta_3+\alpha_2,\alpha_2,\alpha_1\}$ in type $F_4$ at $l_i=4$.\vspace{2mm}
\item $\Delta_{l_i}=\{\beta_2+\alpha_1,\alpha_1\}$ in type $G_2$ at $l_i=3$.
\end{itemize}
Finally, for a factor of type $G_2$ at $l_i=2$ we take $\Delta_{l_i}=\{2\alpha+\beta,\beta,\alpha\}$.  This covers all cases, and we now have the explicit subset of positive roots
\[
\Delta_l:=\coprod_{i\in \pi_0(\Delta)}\Delta_{l_i}.
\]
\par

According to Propositions \ref{prop:863} and \ref{prop:865}, the roots $\alpha\in \Delta_l$ label the skew-primitive generators $E_\alpha$ and $F_\alpha$ in $\mfk{v}_q$, for arbitrary $G$ at an arbitrary torsion parameter $q$.  One can check that $\Delta_l$ also provides a basis for the subsystem
\begin{equation}\label{eq:1028}
\{\gamma\in \Phi:l_\gamma>1\}\ \subseteq\ \Phi,
\end{equation}
whenever $G$ has no factor of type $G_2$ at $l_i=2$.  In all cases, the subsystem \eqref{eq:1028} labels the root vectors $E_\gamma$ and $F_\gamma$ which appear in $\mfk{v}_q$.

\subsection{Small quantum algebras below the lacing number}
\label{sect:smallest2}

Having taken proper account of the above details, we can now proceed exactly as in Section \ref{sect:sqa}.
\par

Take
\[
A=(X/2\cdot lQ)^\vee,\ \ \bar{A}=(X/X^\ast)^\vee,
\]
and consider the characters $K_\nu:X\to k^\ast$, $K_\nu(\lambda)=q(\lambda,\nu)$.  These characters vanish on $2\cdot lQ$, and so define elements in $A$.
\par

We note that $\mfk{v}_q$ is a Hopf subalgebra in $U_q$ and we have the modification of the grouplikes realized by the algebra
\[
v_q:=\frac{\mfk{v}_q\rtimes A}{(\msf{K}_\nu\cdot K_\nu^{-1}:\nu\in Q)}.
\]
The algebra $v_q$ inherits a unique Hopf structure so that the characters $A$ are grouplike and the map $\mfk{v}_q\to v_q$ is a map of Hopf algebras.  We additionally have a triangular decomposition for $v_q$ which is inherited from that of $\mfk{v}_q$.

\begin{definition}
Define $\bar{u}_q=\bar{u}(G_q)$ to be the subalgebra in $v_q$ which is generated by the $\mbf{E}_\alpha=K_\alpha E_\alpha$ and $F_\alpha$, for $\alpha\in \Delta_l$, and the subgroup of characters $\bar{A}\subseteq A$.
\end{definition}

We again refer to $\bar{u}_q$ as the \emph{smallest quantum algebra} for $G$ at $q$.  We have the triangular decomposition, distinguished basis, and dimension calculation
\[
\dim(\bar{u}_q)=[X:X^\ast]\cdot (\prod_{\gamma\in \Phi^+}l_\gamma)^2,
\]
just as in Section \ref{sect:sqa}.

\subsection{Realization via coinvariants}

We describe an alternate means of locating the subalgebra $\bar{u}_q$ in $v_q$.  Consider the group of characters $\Psi=(X^\ast/2\cdot lQ)^\vee$ along with the natural quotient map $A\to \Psi$, which is given by restriction.
\par

The composition
\[
\Rep(G^\ast_{\varepsilon})\to \Rep(G_q)\overset{\res}\longrightarrow \Rep(v_q)
\]
provides a tensor functor from $\Rep(G^\ast_{\varepsilon})$ to $\Rep(v_q)$, and each $G^\ast_{\varepsilon}$-representation decomposes into a sum of $1$-dimensional representations over $v_q$.  These representations are labeled by elements in the quotient $X^\ast/2\cdot lQ$ and define a map of Hopf algebras $\pi:v_q\to k\Psi$.  The map $\pi$ annihilates all of the generators $E_\alpha$ and $F_\alpha$, and lifts the group map $A\to \Psi$.
\par

Via Nichols-Zoeller freeness \cite{nicholszoeller89} the subalgebra of left $k\Psi$-coinvariants in $v_q$ is seen to be of dimension
\[
\dim({^{k\Psi}v_q})=\frac{\dim(v_q)}{|\Psi|}=[X:X^\ast]\cdot (\prod_{\gamma\in \Phi^+}l_\gamma)^2.
\]
Furthermore, by checking on the generators, one observes an algebra inclusion $\bar{u}_q\subseteq {^{k\Psi}v_q}$.  Via the triangular decomposition for $v_q$ one sees that $\bar{u}_q$ has dimension greater than or equal to $[X:X^\ast]\cdot (\prod_{\gamma\in \Phi^+}l_\gamma)^2$, which forces a calculation
\[
\dim(\bar{u}_q)=[X:X^\ast]\cdot (\prod_{\gamma\in \Phi^+}l_\gamma)^2,\ \ \text{an equality of algebras}\ \ \bar{u}_q= {^{k\Psi}v_q},
\]
and verifies the triangular decomposition for $\bar{u}_q$ as well.

\subsection{Tensor structures}
\label{sect:res_ten2}

As before, the algebra $\bar{u}_q$ sits in $v_q$ as a right coideal subalgebra in $v_q$.  This follows from the computations
\begin{equation}\label{eq:856}
\Delta(\mathbf{E}_\alpha)=\mathbf{E}_\alpha\ot K_\alpha+K^2_\alpha\ot \mathbf{E}_\alpha\ \ \text{and}\ \ \Delta(F_\alpha)=F_\alpha\ot K^{-1}_\alpha+1\ot F_\alpha,
\end{equation}
for $\alpha\in \Delta_l$, and the fact that all $K_\alpha^2|_{X^\ast}= 1$.
\par

Via the coaction \eqref{eq:856} we obtain a natural right action of $\Rep(G_q)$ on $\Rep(\bar{u}_q)$, and the restriction functor
\[
\res^-_q:\Rep(G_q)\to \Rep(\bar{u}_q)
\]
becomes a map of $\Rep(G_q)$-module categories.  So, $\Rep(\bar{u}_q)$ again carries a natural pointed module category structure over $\Rep(G_q)$.  The details here are exactly as in Section \ref{sect:res_ten}.

\subsection{Braid group actions}

By construction, the action of the braid group $B_{\mfk{g}}$ on $U_q$ restricts to an action of $B_{\mfk{g}}$ on the Hopf subalgebra $\mfk{v}_q$.  This action extends uniquely to an action on $v_q$ under which the map $\mfk{v}_q\to v_q$ is $B_{\mfk{g}}$-equivariant, and under which $B_{\mfk{g}}$ acts on the grouplikes $\xi\in A$ as $T_\alpha\cdot \xi=\xi(\sigma^{-1}_\alpha-)$.  We claim that the subalgebra $\bar{u}_q\subseteq v_q$ is stable under this braid group action.

\begin{proposition}
The algebra $\bar{u}_q$ admits a natural action of the braid group $B_{\mfk{g}}$ under which the restriction functor $\res^-_q:\Rep(G_q)\to \Rep(\bar{u}_q)$ is $B_\mfk{g}$-equivariant.
\end{proposition}

For the proof one employs the identification of $\bar{u}_q$ with the $k\Psi$-coinvariants in $v_q$, and leverages $B_{\mfk{g}}$-equivariance of the map $\pi:v_q\to k\Psi$.  We leave the details to the interested reader.

\section{Normality of restriction functors}
\label{sect:normal}

In this section we show that the restriction functor $\res^-_q:\Rep(G_q)\to \Rep(\bar{u}_q)$  is normal.  Here $\Rep(\bar{u}_q)$ is considered along with its natural structure as a pointed module category over $\Rep(G_q)$ (see Sections \ref{sect:res_ten} and \ref{sect:res_ten2}).  We show additionally that the kernel of this functor is precisely the image of quantum Frobenius in $\Rep(G_q)$.
\par

At certain parameters such normality is well-known.  See for example \cite{deconcinilyubashenko94} \cite[\S\ 2.3]{davydovetingofnikshych18} for the case of simple $G$ at an odd order scalar parameter which is coprime to the determinant of the Cartan matrix.  The general odd order case also follows from standard relations which can be found in \cite{lusztig90II,lusztig93}.  One can see \cite[Proposition 1.5]{arkhipovgaitsgory03} for the case of a semisimple simply-connected group $G$ at a quantum parameter which is of ``sufficiently large" even order.\footnote{As far as we can tell, one also assumes that a sufficiently large power of $2$ divides all $l_i$ in \cite{arkhipovgaitsgory03}.}

\subsection{Normality of restriction}
\label{sect:normal.1}

\begin{proposition}[{cf.\ \cite[Problem 7.1]{lentner16}}]\label{prop:norm}
The restriction functor $\res^-_q:\Rep(G_q)\to \Rep(\bar{u}_q)$ is normal.  Furthermore, the kernel of $\opn{res}^-_q$ is equal to the image of quantum Frobenius
\[
\opn{ker}(\res^-_q)=\opn{im}\big(\opn{Fr}:\Rep(G^\ast_{\varepsilon})\to \Rep(G_q)\big).
\]
\end{proposition}

The proposition is obtained as an application of Lemma \ref{lem:norm1} below.  We provide a proof which is contingent on Lemma \ref{lem:norm1}, then return to cover the necessary details.

\begin{proof}
All objects in $\opn{ker}(\res^-_q)$ necessarily have grading in the kernel of the projection $X\to (A)^\vee=X/X^\ast$, and are annihilated by all $E_\alpha$ with $l_\alpha>1$.  So $\opn{ker}(\res^-_q)\subseteq \opn{im}(\opn{Fr})$ by Lemma \ref{lem:696} and Theorem \ref{thm:LQF}.  We refer again to Lemma \ref{lem:696} to see that $\opn{im}(\opn{Fr})\subseteq \opn{ker}(\res^-_q)$, and hence that the kernel is precisely the image of quantum Frobenius.
\par

We need to show now that the $\bar{u}_q$-invariants $V^{\bar{u}_q}$ in any $G_q$-representation $V$ form a $G_q$-subrepresentation.  We have explicitly
\[
V^{\bar{u}_q}=\left\{\begin{array}{c}
\text{The sum of all homogeneous $v\in V$ with}\\
\deg(v)\in X^\ast\text{ and }\mathbf{E}_\alpha\cdot v=F_\alpha\cdot v=0\text{ at all }\alpha\in \Delta_l
\end{array}
\right\}\vspace{2mm}
\]
\[
\hspace{1cm}=\left\{\begin{array}{c}
\text{The sum of all homogeneous $v\in V$ with}\\
\deg(v)\in X^\ast\text{ and }E_\alpha\cdot v=F_\alpha\cdot v=0\text{ at all }\alpha\in \Delta_l
\end{array}
\right\}.
\]
\par

By definition, this subspace $V^{\bar{u}_q}$ is stable under the actions of $E_\alpha$ and $F_\alpha$, for $\alpha\in \Delta_l$, and we need only show it is stable under the actions of the divided powers $E_\beta^{(l_\beta)}$ and $F_\beta^{(l_\beta)}$ for simple $\beta$.  Rather, we need to show that, for any $v\in V^{\bar{u}_q}$, the vectors
\[
E^{(l_\beta)}_\beta\cdot v\ \ \text{and}\ \ F_\beta^{(l_\beta)}\cdot v
\]
are still annihilated the actions of the elements $E_\alpha,\ F_\alpha$, and $(\xi-1)$ with $\xi\in \bar{A}$.  However, this is precisely what Lemma \ref{lem:norm1} tells us.
\end{proof}

We now cover the requisite details.  In our arguments below we consider the topological Hopf algebra
\[
\hat{U}_q:=\varprojlim_V \dot{U}_q/\opn{Ann}(V),
\]
where the limit is over all finite-dimensional $G_q$-representations $V$.  Here we employ the partial ordering where $V\leq W$ whenever $W$ admits an injective map of representations $V\to W$.  Alternatively, $\hat{U}_q$ is constructed as the endomorphism algebra of the forgetful functor
\[
\hat{U}_q=\opn{End}\big(forget:\Rep G_q\to \opn{Vect}\big)
\]
\cite[\S 2.6]{negron21}.  This algebra is generated topolgically by the divided powers of the simple root vectors $E^{(n)}_\beta$ and $F^{(n)}_\beta$, and the characters $X^\vee=\opn{Hom}_{\opn{Grp}}(X,k^{\times})$.
\par

By Tannakian reconstruction, the restriction functor $\Rep(G_q)\to \Rep(\bar{u}_q)$ defines an algebra map $\bar{u}_q\to \hat{U}_q$, which explicitly sends the $\mbf{E}_\gamma$ and $F_\gamma$ in $\bar{u}_q$ to the corresponding vectors in $\hat{U}_q$ and the characters $\xi\in \bar{A}$ to their corresponding functions on $X$.

\begin{remark}
From a simple-minded perspective, $\hat{U}_q$ is just a version of the quantum enveloping algebra which is unital--unlike the modified algebra $\dot{U}_q$--and which contains all naturally occurring toral characters--unlike the usual quantum enveloping algebra $U_q$.  In $\hat{U}_q$ we can directly compare elements from all of our favorite algebras; $U_q$, $\mfk{v}_q$, and $\bar{u}_q$.
\end{remark}

\begin{lemma}\label{lem:norm1}
Let $\mfk{m}\subseteq \bar{u}_q$ be the kernel of augmentation $\bar{u}_q\to k$.  For all $\xi\in \bar{A}$, and all positive roots $\gamma$, we have
\begin{equation}\label{eq:1018}
(\xi-1)E_\gamma^{(l_\gamma)}=E_\gamma^{(l_\gamma)}(\xi-1)\ \ \text{and}\ \ (\xi-1)F_\gamma^{(l_\gamma)}=F_\gamma^{(l_\gamma)}(\xi-1).
\end{equation}
Additionally, for all $\alpha\in \Delta_l$ and all $\beta\in \Delta$ we have
\begin{equation}\label{eq:1006}
E_{\alpha}E^{(l_\beta)}_{\beta},\ F_\alpha F^{(l_\beta)}_\beta,\ E_{\alpha}F^{(l_\beta)}_{\beta},\ F_\alpha E^{(l_\beta)}_\beta \ \in\  \hat{U}_q\mfk{m}.
\end{equation}
\end{lemma}

The result essentially follows by computations of Lusztig \cite{lusztig90II,lusztig93} and Lentner \cite{lentner16}.

\begin{proof}
It suffices to prove the result when $G$ is almost-simple.  Hence we take $l=l_i$ at the unique index $i\in \pi_0(\Delta)$.  When $l=1$ the group $\bar{A}$ is trivial and $\Delta_l$ is empty, so that $\bar{u}_q=k$ and there is nothing to check.  We therefore suppose that $l>1$.
\par

Before beginning, let us note that we have containments
\[
\hat{U}_qE_\gamma\ \subseteq\ \hat{U}_q\mfk{m}\ \ \text{and}\ \ \hat{U}_qF_\gamma \subseteq\ \hat{U}_q\mfk{m}
\]
whenever $\gamma$ is a positive root with $l_\gamma>1$.  We take these inclusions for granted, and use them throughout the proof.
\par

Each characters $\xi$ in $(X/X^\ast)^\vee=\bar{A}$ has $\xi|_{l Q}= 1$, by Lemma \ref{lem:Q_v_X}.  Hence $\xi(l_\gamma\gamma)=1$ at each root $\gamma$, and the commutativity relations \eqref{eq:1018} follow.
\par

As for the relations \eqref{eq:1006}, it suffices to establish inclusions
\begin{equation}\label{eq:1045}
\left.\begin{array}{r}
{[E_\beta^{(l_\beta)},E_\alpha]},\ \ {[F_\beta^{(l_\beta)},F_\alpha]}\ \\
{[E_\beta^{(l_\beta)},F_\alpha]},\ \ {[F_\beta^{(l_\beta)},E_\alpha]}\ 
\end{array}\right\}\ \in\ \hat{U}_q\mfk{m}
\end{equation}
at arbitrary $\alpha\in \Delta_l$ and $\beta\in \Delta$.  When both $\alpha,\beta\in \Delta$, and $l_\alpha,l_\beta >1$, one employs the higher Serre relations \cite[Proposition 35.3.1]{lusztig93} to see that the first two commutators lie in $\hat{U}_q\mfk{m}$.  Similarly, in this case the relations of \cite[\S 6.5 (a2), (a6)]{lusztig90II} tell us that the latter two commutators lie in $\hat{U}_q\mfk{m}$ as well.  We therefore have both of the desired containments \eqref{eq:1018} and \eqref{eq:1006}, provided $l$ is greater than the lacing number for $G$.
\par

The cases where $l$ is less than or equal to the lacing number for $G$ are somewhat involved, and are delayed until Appendix \ref{sect:norm_app}.
\end{proof}

\subsection{Normality and Hopf structures}

We note that the fiber functor for $\Rep(G_q)$ endows $\opn{Vect}$ with the structure of a pointed $\Rep(G_q)$-module category, and that the forgetful functor $\Rep(\bar{u}_q)\to \opn{Vect}$ is naturally a map of pointed module categories over $\Rep(G_q)$.  It follows that the restriction functor $\res^-_q:\Rep(G_q)\to \Rep(\bar{u}_q)$ determines, and is determined by, a map of $\O(G_q)$-module coalgebras $\pi:\O(G_q)\to \bar{u}_q^\ast$ \cite[Definition 1.1]{doi92}, via Tannakian reconstruction \cite[Theorem 2.2.8]{schauenburg92}.
\par

Similarly, the quantum Frobenius functor determines an embedding of Hopf algebras																		 $fr^\ast:\O(G^\ast_{\varepsilon})\to \O(G_q)$ \cite[Lemma 2.2.13]{schauenburg92}.  The following codifies normality of the restriction functor in algebraic terms.

\begin{lemma}\label{lem:coinv}
Consider a Hopf algebra $A$, and let $\pi:A\to C$ be a map of $A$-module coalgebras for which the restriction functor $\res_\pi:\opn{Corep}(A)\to \opn{Corep}(C)$ is normal.  Let $B$ be the (unique) Hopf subalgebra in $A$ whose corepresentations are identified with the kernel of $\opn{res}_\pi$.  Then the left and right $C$-coinvariants in $A$ agree, and ${^C A}=A^C=B$.
\end{lemma}

\begin{proof}
The fact that $B$ is in the kernel of $\res_\pi$ implies $B\subseteq A^C$, by definition.  The left counit axiom for $B$ then gives $\pi|_B=\epsilon_B$.  Explicitly, for $b\in B$  we have
\[
\epsilon(b)=\epsilon(b_1)\pi(b_2)=\pi(\epsilon(b_1)b_1)=\pi(b).
\]
Since $B$ is a subcoalgebra in $A$, this implies $B\subseteq {^CA}$ as well.
\par

For the opposite inclusion $A^C\subseteq B$, normality of restriction tells us that $A^C$ lies in $\opn{Corep}(B)$, so that $\Delta(A^C)\subseteq A\ot B$.  We apply the counit on the left to find $A^C\subseteq B$, and hence $A^C=B$.
\par

The final inclusion ${^C A}\subseteq B$ follows from the $A$-module coalgebra structure on $C$.  Indeed, for any $a\in {^CA}$ we have $\pi(a_1)\ot a_2=1\ot a\in C\ot A$, and hence
\[
\pi(S a_1)\ot a_2=\big(1\cdot Sa_1\big)\ot a_2=\big(\pi(a_1)\cdot Sa_2\big)\ot a_3=1\ot a.
\]
Apply $id_C\ot S$ to the above expression to see that $S(a)\in A^C=B$ whenever $a\in {^CA}$.  This gives $a=S^{-1}S(a)\in S^{-1}(B)=B$, and verifies the equality ${^CA}=B$.
\end{proof}

We apply Lemma \ref{lem:coinv} to the quantum group setting, as described above.

\begin{corollary}\label{cor:coinv}
For arbitrary $G$ and $q$,
\[
{^{\bar{u}^\ast_q}\O(G_q)}=\O(G_q)^{\bar{u}^\ast_q}=\O(G^\ast_{\varepsilon}).
\]
\end{corollary}

\subsection{Normality of restriction for the quantum Borel}

We consider the category $\Rep(B_q)$ of representations for the positive quantum Borel.  These are integrable $U_q^+$-representations which are equipped with a compatible grading by the character lattice $X$ for $G$.  As with the full quantum group, we have the tensor embedding
\[
\opn{Fr}^{\geq 0}:\Rep(B_\varepsilon^\ast)\to\Rep(B_q)
\]
which is obtained by restricting along the algebra map $fr|_{U^+_q}:U^+_q\to (U^\ast_\varepsilon)^+$ \cite{lusztig93,lentner16}.
\par

We also have the non-negative subalgebra $\bar{u}^{\geq 0}_q$ in $\bar{u}_q$, which sits as a subcomodule subcoalgebra in $\bar{u}_q$.  Hence we obtain a natural action of $\Rep(B_q)$ on the category of representations $\Rep(\bar{u}^{\geq 0}_q)$, and the restriction functor
\[
\Rep(B_q)\to \Rep(\bar{u}^{\geq 0}_q)
\]
endows $\Rep(\bar{u}_q^{\geq 0})$ with a pointed module category structure over $\Rep(B_q)$.
\par

One can use Lentner's normality result \cite[Lemma 6.6]{lentner16}, or \cite[Proposition 35.3.1]{lusztig93}, to obtain the following non-negative analog of Proposition \ref{prop:norm}.

\begin{proposition}[{\cite{lusztig93,lentner16}}]\label{prop:B_norm}
The restriction functor $\Rep(B_q)\to \Rep(\bar{u}^{\geq 0}_q)$ is normal, and has kernel equal to the image of $\Rep(B^\ast_\varepsilon)$ in $\Rep(B_q)$ under quantum Frobenius.
\end{proposition}

\section{Restriction and simple $G_q$-representations}

Before turning to a discussion of the Steinberg representation, we record some helpful results which describe simple $G_q$-representations and their behaviors under restriction.  For convenience of presentation,
\begin{quote}
\emph{We suppose that $G$ is simply-connected throughout this section}.
\end{quote}

\begin{remark}
One recovers an analysis of the simples for non-simply-connected $G$ by noting that $\Rep(G_q)$ sits in its simply-connected form $\Rep(G^{sc}_q)$ as the full subcategory of representations whose $P$-grading is supported on the character lattice $X$ for $G$.
\end{remark}

\subsection{Generic information}

The following is well-known at most parameters $q$.

\begin{lemma}[{\cite[Lemma 1.13]{andersenpolowen91}}]\label{lem:weights}
Consider a $G_q$-representation $V$, and let $\Pi_l(V)$ denote the subset of weights $\mu\in P$ with $V_\mu$ non-vanishing.  Then $\Pi_l(V)$ is stable under the action of the Weyl group on $P$.
\end{lemma}

\begin{proof}
It suffices to consider the case where $V$ is finite-dimensional.  Let $v\in V$ be nonzero and of weight $\mu$.  Fix simple $\alpha$ and take $m=-(\alpha,\mu)/d_\alpha$.  We claim that if $m\geq 0$ then $E^{(m)}_\alpha\cdot v\neq 0$, and if $m\leq 0$ then $F^{(m)}_\alpha\cdot v\neq 0$.  Supposing this claim is valid, we find that $V_{\sigma_\alpha(\mu)}=V_{\mu+m\alpha}$ is non-zero, that $\Pi_l(V)$ is stable under the actions of the simple reflections in $\msc{W}$, and hence that $\Pi_l(V)$ is stable under the action of the Weyl group.
\par

To establish the claim, we reduce to the case of $G=\opn{SL}(2)$ at a parameter $\zeta$ by restricting along the root subgroup for $\alpha$.  Now, by induction on the length, we reduce further to the case where $V$ is a simple $\opn{SL}(2)_\zeta$-representation.  Here all the simples are of the form $\opn{Fr}(L_1)\ot L_2$ where $L_1$ is a representation for classical $\SL(2)$ and $L_2$ has highest weight $<\opn{ord} \zeta(\alpha,\alpha)$, and the claim follows by a direct inspection.
\end{proof}

For a given simple representation $L(\lambda)$ we take $\Pi_l(\lambda)=\Pi_l(L(\lambda))$.  A classical fact tells us that the minimal element in this set is $w_0\lambda$ \cite[Theorem 4.12]{lusztig88} \cite[Proposition 21.3, Exercise 21.4(6)]{humphreys12}.  We therefore obtain the following.

\begin{lemma}\label{lem:lowest_wt}
For $\lambda\in P^+$, $L(\lambda)$ has lowest weight $w_0\lambda$, and $L(\lambda)^\ast=L(-w_0\lambda)$.
\end{lemma}

\subsection{Restrictions of simple $G_q$-representations}

Below we take
\[
P_l=\left\{\lambda\in P:0\leq \frac{(\alpha,\lambda)}{d_\alpha}<l_\alpha\ \forall\ \alpha\in \Delta\right\}.
\]
In \cite{parshallwang91} such weights are called \emph{($l$-)restricted}.  These are weights of the form $\lambda=\sum_{\alpha}c_\alpha \omega_\alpha$, where the $\omega_\alpha$ are the fundamental weights associated to simple roots $\alpha$, and the coefficients satisfy $0\leq c_\alpha <l_\alpha$.  We note that the class of simples $L(\lambda)$ of restricted highest weight is closed under duality.  This just follows from the fact that $w_0\Delta=-\Delta$.

\begin{lemma}[{\cite[Lemma 9.3.2]{parshallwang91}}]\label{lem:1366}
For any $\lambda\in P_l$, the simple representation $L(\lambda)$ satisfies
\[
L(\lambda)^{\bar{u}^+_q}=L(\lambda)^{\mfk{v}^+_q}=L(\lambda)_{\lambda}.
\]
\end{lemma}

To be clear, by the $\bar{u}^+_q$ or $\mfk{v}^+_q$-invariants here we mean the collections of vectors which are annihilated by all $\mbf{E}_\alpha$ or $E_\alpha$, for $\alpha\in \Delta_l$, respectively.

\begin{proof}
The equality $L(\lambda)^{\bar{u}^+_q}=L(\lambda)^{\mfk{v}^+_q}$ is immediate, since all $K_\gamma$ act as units on $L(\lambda)$.  We restrict to $B_q$ and, for any $\mu\in P$, consider the product $L(\lambda)\ot k_\mu$.  It suffices to show that the invariants $(L(\lambda)\ot k_\mu)^{\bar{u}^{\geq 0}_q}$ lie in the span of the highest weight vector, at arbitrary $\mu\in P$.  Let $v(\mu)$ denote the (unique-up-to-scaling) highest weight vector in $L(\lambda)\ot k_\mu$, which is just the highest weight vector in $L(\lambda)$ times any nonzero vector in $k_\mu$.

By the normality result of Proposition \ref{prop:B_norm}, the invariants $(L(\lambda)\ot k_\mu)^{\bar{u}^{\geq 0}_q}$ are identified with a $B^\ast_\varepsilon$-representation via quantum Frobenius, and the socle of this representation over $B^\ast_\varepsilon$ lies in the socle of $L(\lambda)\ot k_\mu$ over $B_q$, i.e.\ in the highest weight space.  It follows that $v(\mu)$ is contained in the $B^\ast_\varepsilon$-subrepresentation generated by any nonzero vector $w$ in $(L(\lambda)\ot k_\mu)^{\bar{u}_q^{\geq 0}}$.
\par

From the expression $\lambda=\sum_{\alpha}c_\alpha\omega_\alpha$, with each $c_\alpha<l_\alpha$ by hypothesis, we find at each simple root $\alpha$ that
\[
\sigma_\alpha(\lambda-l_\alpha\alpha)=\lambda-c_\alpha\alpha+l_\alpha\alpha=\lambda+(l_\alpha-c_\alpha)\alpha >\lambda.
\]
In particular $\lambda-l_\alpha\alpha\notin \Pi_l(\lambda)$.  It follows that any nonzero homogeneous $w$ in $(L(\lambda)\ot k_\mu)^{\bar{u}_q^{\geq 0}}$ which is of degree $<\lambda+\mu$ generates a subrepresentation which does \emph{not} contain $v(\mu)$.  By our previous observations about $B^\ast_{\varepsilon}$-subrepresentations in $(L(\lambda)\ot k_\mu)^{\bar{u}_q^{\geq 0}}$, no such nonzero vector exists.  This forces
\[
(L(\lambda)\ot k_\mu)^{\bar{u}_q^{\geq 0}}\ \subseteq\ k\cdot v(\mu).
\]
\end{proof}

\begin{proposition}[{\cite[Proposition 9.3.4]{parshallwang91}}]\label{prop:res_simples}
For each $\lambda\in P_l$, the associated simple $L(\lambda)$ in $\Rep(G_q)$ restricts to a simple object in $\Rep(\bar{u}_q)$.
\end{proposition}

\begin{proof}
One employs Lemma \ref{lem:1366} and proceeds exactly as in \cite{parshallwang91}.
\end{proof}

By the usual analysis with baby Verma modules, one sees that simple representations in $\Rep(\bar{u}_q)$ are characterized by their highest weights.  These highest weights are, by definition, characters for the group $\bar{A}=(P/P^\ast)^\vee$, i.e.\ elements in the quotient $P/P^\ast$.  Now, the projection $P\to P/P^\ast$ restricts to provide a bijection $P_l\overset{\sim}\to P/P^\ast$, so that
\[
|\opn{Irrep}(\bar{u}_q)|=|P/P^\ast|=|P_l|.
\]
Hence Proposition \ref{prop:res_simples} implies the following.

\begin{corollary}\label{cor:simples}
For any simple representation $L(\bar{\lambda})$ in $\Rep(\bar{u}_q)$ there is a simple representation $L(\lambda)$ in $\Rep(G_q)$ for which
\[
\res^-_q(L(\lambda))=L(\bar{\lambda}).
\]
\end{corollary}

One sees from Corollary \ref{cor:lus_st} below that this simple lift $L(\lambda)$ of $L(\bar{\lambda})$ is unique.

\subsection{Steinberg decomposition}

We use Lemma \ref{lem:1366} and Proposition \ref{prop:res_simples} to recover Steinberg's tensor product theorem at arbitrary $q$.

\begin{theorem}[{\cite[Theorem 7.4]{lusztig89}}]\label{thm:lus_st}
Suppose that $G$ is simply-connected.  For any $\lambda\in P^+$ we have a unique decomposition $\lambda=\lambda_0+\lambda_1$ with $\lambda_1\in P^\ast$ and $\lambda_0\in P_l$, and a corresponding isomorphism of $G_q$-representations
\[
\opn{Fr}(L(\lambda_1))\ot L(\lambda_0)= L(\lambda_1)\ot L(\lambda_0)\overset{\sim}\to L(\lambda).
\]
\end{theorem}

\begin{proof}
We check highest weight vectors in the product $L(\lambda_1)\ot L(\lambda_0)$.  Since $L(\lambda_1)$ is in the kernel of the restriction functor $\Rep(B_q)\to \Rep(\mfk{v}^+_q)$ we have
\[
(L(\lambda_1)\ot L(\lambda_0))^{\mfk{v}^+_q}=L(\lambda_1)\ot L(\lambda_0)^{\mfk{v}^+_q}=L(\lambda_1)\ot k_{\lambda_0}.
\]
The final equality here follows by Proposition \ref{prop:res_simples}.  We note that $ L(\lambda_1)\ot k_{\lambda_0}$ has $1$-dimensional socle over $B_q$ to see now that $L(\lambda_1)\ot L(\lambda_0)$ has a $1$-dimensional socle over $B_q$.  This $B_q$-socle is spanned by the unique-up-to-scaling element of weight $\lambda$.  This calculates the socle over $G_q$ as $\opn{soc}(L(\lambda_1)\ot L(\lambda_0))=L(\lambda)$.
\par

We note at this point that $L(\lambda_1)\ot L(\lambda_0)$ and $L(\lambda)$ both have lowest weights $w_0\lambda$.  Take duals and apply a similar analysis to see that $L(\lambda_0)^\ast\ot L(\lambda_1)^\ast$ has socle $L(-w_0\lambda)$.  We therefore observe a projection
\[
p:L(\lambda_1)\ot L(\lambda_0)\to L(-w_0\lambda)^\ast=L(\lambda).
\]
The kernel of this map must be trivial, since $L(\lambda_1)\ot L(\lambda_0)$ has no additional highest weight vectors, outside of the one of weight $\lambda$.  Hence $p$ is injective.  Surjectivity follows by simplicity of $L(\lambda)$, so that we observe the claimed isomorphism.
\end{proof}

From Proposition \ref{prop:res_simples} we know that $\res^-_q(L(\lambda))$ is simple in $\Rep(\bar{u}_q)$ whenever $\lambda\in P_l$.  So the Steinberg decomposition can alternately be expressed in terms of the restriction functor $\res^-_q$.

\begin{corollary}\label{cor:lus_st}
Suppose that $G$ is simply-connected.  Then for any $\lambda\in P^+$ the simple representation $L(\lambda)$ admits a unique decomposition
\[
L(\lambda)=L(\lambda'')\ot L(\lambda')
\]
where $L(\lambda')$ restricts to a simple representation in $\Rep(\bar{u}_q)$, and $L(\lambda'')$ is in the kernel of $\res^-_q$.
\end{corollary}

\section{The Steinberg representation}
\label{sect:steinborg}

We consider the dominant weight $\rho_l:=\frac{1}{2}\sum_{\gamma\in \Phi^+}(l_\gamma-1)\gamma$ and its corresponding simple representation
\[
\opn{St}:=L(\rho_l)
\]
in $\Rep(G_q)$, for simply-connected $G$.  The representation $\opn{St}$ is known as the \emph{Steinberg representation}, and it has proved extraordinarily useful in studies of both modular and quantum representation theory.  See for example \cite{jantzen03,parshallwang91,andersenpolowen91,andersenpolowen92,andersen03}.
\par

Our next goal is to show that the Steinberg representation is both projective and injective in $\Rep(G_q)$, and that it restricts to a projective and injective object in $\Rep(\bar{u}_q)$  (cf.\ \cite{andersen03,parshallwang91,andersenpolowen91,andersenpolowen92,andersen92}).  However, these points are more easily argued in a locale which lies between these two categories.
\par

In this section we analyze the behaviors of the Steinberg representation in the category of \emph{torally extended} small quantum group representations $\Rep(\dot{\mfk{u}}_q)$.   We then return to a discussion of $G_q$-representations and $\bar{u}_q$-representations in Sections \ref{sect:steinborg_sc} and \ref{sect:steinborg_nsc}.  Our approach follows those of Parshall-Wang \cite{parshallwang91} and Andersen-Polo-Wen \cite{andersenpolowen91}, to an extent, though we employ basic results for Hopf algebras \cite{larsonsweedler69} in order to avoid intricate arguments with induction.
\par

We note that $\rho_l$ does not lie in the character lattice for $G$, at a general non-simply-connected group.  So we restrict our attention to the simply-connected setting.
\begin{center}
\emph{Throughout the section we assume $G$ is simply-connected!}
\end{center}

\subsection{Torally extended small quantum group representations}
\label{sect:dotuq}

We consider the tensor category $\Rep(\dot{\mfk{u}}_q)$ of ``torally extended" small quantum group representations.  Specifically we consider
\begin{equation}\label{eq:1391}
\Rep(\dot{\mfk{u}}_q)=\left\{\begin{array}{c}
\text{The category of $\mfk{v}_q$-modules with a}\\
\text{compatible grading by the character lattice }X
\end{array}\right\}.
\end{equation}
Since we've restricted ourselves to the simply-connected setting, this character lattice is just the weight lattice $X=P$.  By a compatible $P$-grading we mean that $E_\gamma$ and $F_\gamma$ shift degrees on a $\dot{\mfk{u}}_q$-representation $V$ by $\gamma$ and $-\gamma$ respectively, and that the $\msf{K}_\nu$ act on each homogeneous subspace $V_\lambda\subseteq V$ as the corresponding scalar $q(\nu,\lambda)$.
\par

Our notation \eqref{eq:1391} is taken from \cite[\S\ 31.1 \& 36.2]{lusztig93}, where $\Rep(\dot{\mfk{u}}_q)$ is constructed as the category of ``unital" representations for a non-unital algebra $\dot{\mfk{u}}_q$ with generators $E_\alpha 1_\lambda$ and $F_\alpha 1_\lambda$, for $\alpha\in \Delta_l$ and $\lambda\in P$, which satisfy
\[
E_\alpha 1_\lambda\cdot F_\beta 1_\mu=\delta_{\mu,\lambda+\beta}E_\alpha F_\beta 1_\mu\ \ \text{and}\ \ F_\beta 1_\mu\cdot E_\alpha 1_\lambda=\delta_{\lambda,\mu-\alpha}F_\beta E_\alpha 1_\lambda,
\]
as well as relations deduced from the usual quantum Serre relations.

\subsection{The smallest quantum algebra and $\dot{\mfk{u}}_q$-representations}

Each of the generators for $\bar{u}_q$ act in the obvious way on $\dot{\mfk{u}}_q$-representations.  In particular, the elements $\mbf{E}_\gamma$ and $F_\gamma$ act on a $\dot{\mfk{u}}_q$-representation via the corresponding elements in $\mfk{v}_q$, and the characters $\xi\in \bar{A}$ act as the semisimple operators
\[
\xi\cdot v=\xi(\lambda)\cdot v\ \ \text{for homogeneous $v\in V$, }\lambda=\deg(v).
\]
In this way we obtain a natural restriction functor $\Rep(\dot{\mfk{u}}_q)\to \Rep(\bar{u}_q)$ and, as in Section \ref{sect:res_ten}, $\Rep(\bar{u}_q)$ has the natural structure of a pointed module category over $\Rep(\dot{\mfk{u}}_q)$.
\par

For the maximal torus $T^\ast$ in the dual group $G^\ast$, we consider the tensor subcategory $\Rep(T^\ast_\varepsilon)$ of $P^\ast$-graded vector spaces in $\Rep(\dot{\mfk{u}}_q)$.  One checks the following directly.

\begin{lemma}\label{lem:norm_dot}
The restriction functor $\Rep(\dot{\mfk{u}}_q)\to \Rep(\bar{u}_q)$ is normal, and has kernel equal to the subcategory of $P^\ast$-graded vector spaces $\Rep(T^\ast_\varepsilon)\subseteq \Rep(\dot{\mfk{u}}_q)$.
\end{lemma}

We also have a natural restriction functor from $\Rep(G_q)$ to $\Rep(\dot{\mfk{u}}_q)$, which is a tensor functor. This tensor functor factors $\res^-_q$ as
\[
\xymatrix{
	& \Rep(\dot{\mfk{u}}_q)\ar[dr] &\\
\Rep(G_q)\ar[ur]\ar[rr]_{\opn{res}^-_q} & & \Rep(\bar{u}_q).
}
\]
The above diagram and Proposition \ref{prop:res_simples} imply the following.

\begin{lemma}\label{lem:res_simples}
For any restricted weight $\lambda\in P_l$, the associated simple $L(\lambda)$ in $\Rep(G_q)$ restricts to a simple object in $\Rep(\dot{\mfk{u}}_q)$.
\end{lemma}

\subsection{Induction for $\dot{\mfk{u}}_q$-representations}
\label{sect:huh}

We consider the category $\Rep(\dot{\mfk{u}}^{\leq 0}_q)$ of character graded representations for the non-positive subalgebra $\mfk{v}^{\leq 0}_q$ in $\mfk{v}_q$, and the restriction functor $\Rep(\dot{\mfk{u}}_q)\to \Rep(\dot{\mfk{u}}_q^{\leq 0})$.  We have the right adjoint to restriction
\[
\Rind:=\Rind_{\dot{\mfk{u}}_q^{\leq 0}}^{\dot{\mfk{u}}_q}:\Rep(\dot{\mfk{u}}^{\leq 0}_q)\to \Rep(\dot{\mfk{u}}_q),
\]
as well as the left adjoint
\[
\Lind:=\Lind_{\dot{\mfk{u}}^{\leq 0}_q}^{\dot{\mfk{u}}_q}:\Rep(\dot{\mfk{u}}_q^{\leq 0})\to \Rep(\dot{\mfk{u}}_q).
\]
One can write $\Lind(V)$ explicitly as $\Lind(V)=\mfk{v}_q\ot_{\mfk{v}_q^{\leq 0}}V\cong\mfk{v}^+_q\ot V$, where the $E_\gamma$ and $F_\gamma$ act via the left factor and the $P$-grading is given by
\[
\deg(E_{\gamma_1}^{m_1}\cdots E_{\gamma_r}^{m_r}\ot v)=\opn{deg}(v)+\sum_i m_i\gamma_i.
\]

\begin{lemma}\label{lem:1361}
For any finite-dimensional $\dot{\mfk{u}}_q^{\leq 0}$-representation $V$,
\[
\Rind(V)=\Lind(V^\ast)^\ast\ \ \text{and}\ \ \Lind(V)=\Rind(V^\ast)^\ast.
\]
In particular, $\Rind$ sends finite-dimensional $\dot{\mfk{u}}_q^{\leq 0}$-representations to finite-dimensional $\dot{\mfk{u}}_q$-representations.
\end{lemma}

In the proof we employ the pivotal structure $id\overset{\sim}\to {(-)^\ast}^\ast$ on $\Rep(\dot{\mfk{u}}_q)$ which is provided by the action of the character $K_\rho$.

\begin{proof}
We consider the first identification.  Since $\Rep(\dot{\mfk{u}}_q)$ is generated by its finite-dimensional representations, under filtered colimits, it suffices to provide natural isomorphisms
\[
\Hom_{\dot{\mfk{u}}_q}(W,\Lind(V^\ast)^\ast)\cong \Hom_{\dot{\mfk{u}}^{\leq 0}_q}(W,V)
\]
at all finite-dimensional $W$ in $\Rep(\dot{\mfk{u}}_q)$.  Now, for finite-dimensional $V$ and $W$ we have the sequence of natural isomorphisms
\[
\begin{array}{rl}
\Hom_{\dot{\mfk{u}}_q}(W,\Lind(V^\ast)^\ast)&\cong \Hom_{\dot{\mfk{u}}_q}({W^\ast}^\ast,\Lind(V^\ast)^\ast)\\
&\cong\Hom_{\dot{\mfk{u}}_q}(\Lind(V^\ast),W^\ast)\vspace{1mm}\\
&\cong\Hom_{\dot{\mfk{u}}^{\leq 0}_q}(V^\ast,W^\ast)\cong\Hom_{\dot{\mfk{u}}^{\leq 0}_q}(W,V)\\
\end{array}
\]
which realize $\Lind(V^\ast)^\ast$ as a, and hence the, right adjoint to restriction.  As a consequence we see that $\Rind$ preserves finite-dimensional representations, since $\Lind$ preserves finite-dimensional representations.  The identification $\Lind(V)=\Rind(V^\ast)^\ast$ is established similarly.
\end{proof}

For any $\lambda\in P$, with associated $1$-dimensional simple $\dot{\mfk{u}}^{\leq 0}_q$-representation $k_\lambda$, we take
\[
\Rind(\lambda)=\Rind(k_\lambda)\ \ \text{and}\ \ \Lind(\lambda)=\Lind(k_\lambda).
\]
We have
\[
\dim(\Rind(\lambda))=\dim(\Lind(\lambda))=\dim(\mfk{v}_q^+)=\prod_{\gamma\in \Phi^+}l_\gamma
\]
at each $\lambda$.  We note that $\Lind(\lambda)$ has a $1$-dimensional lowest weight space of weight $\lambda$, and $1$-dimensional highest weight space of weight $\lambda+2\rho_l$, where $\rho_l=\frac{1}{2}\sum_{\gamma\in \Phi^+}(l_\gamma-1)\gamma$.  Lemma \ref{lem:1361} now tells us that $\Rind(\lambda)$ has $1$-dimensional highest and lowest weights spaces as well, of weights $\lambda$ and $\lambda-2\rho_l$ respectively.

\begin{lemma}\label{lem:1370}
For any $\lambda\in P$, $\Rind(\lambda)$ has a unique highest weight vector which is of weight $\lambda$.  This highest weight vector generates a simple subrepresentation in $\Rind(\lambda)$, and identifies the socle as
\[
\opn{soc}(\Rind(\lambda))=\text{\rm the unique simple $\dot{\mfk{u}}_q$-rep of highest weight }\lambda.
\]
\end{lemma}

\begin{proof}
Let $L'(\mu)$ denote the simple $\dot{\mfk{u}}_q$-representation of highest weight $\mu$.  Highest weight vectors in $\Rind(\lambda)$, of weight $\mu\in P$ say, are identified with maps $M(\mu)\to \Rind(\lambda)$ from the corresponding baby Verma module $M(\mu)=\Lind_{\dot{\mfk{u}}_q^{\geq 0}}^{\dot{\mfk{u}}_q}(\mu)$.  We have by adjunction
\[
\Hom_{\dot{\mfk{u}}_q}(M(\mu),\Rind(\lambda))\cong\Hom_{\dot{\mfk{u}}^{\leq 0}_q}(M(\mu),k_\lambda).
\]
Since $M(\mu)$ is cyclically generated over $\dot{\mfk{u}}^{\leq 0}_q$ by its highest weight space $M(\mu)_{\mu}\cong k_\mu$ we restrict along the projection $M(\mu)\to k_\mu$ to get
\[
\Hom_{\dot{\mfk{u}}_q^{\leq 0}}(M(\mu),k_\lambda)\overset{\sim}\leftarrow\Hom_{\dot{\mfk{u}}_q^{\leq 0}}(k_\mu,k_\lambda)=\left\{
\begin{array}{l}
k\ \ \text{if }\mu=\lambda\\
0\ \ \text{otherwise},\end{array}\right.
\]
and also
\[
\Hom_{\dot{\mfk{u}}^{\leq 0}_q}(k_\mu,k_\lambda)\overset{\sim}\to \Hom_{\dot{\mfk{u}}^{\leq 0}_q}(L'(\mu),k_\lambda)\cong\Hom_{\dot{\mfk{u}}_q}(L'(\mu),\Rind(\lambda)).
\]
The first expression tells us that $\Rind(\lambda)$ admits a highest weight vector of weight $\mu$ if and only if $\lambda=\mu$.  The second expression says that this vector defines a map $L'(\lambda)\to \Rind(\lambda)$ which is unique up to scaling.
\end{proof}

We similarly calculate the cosocle for $\Rind(\lambda)$.

\begin{lemma}\label{lem:cosoc}
At any weight $\lambda$ we have
\[
\opn{cosoc}(\Rind(\lambda))=\text{\rm the unique simple $\dot{\mfk{u}}_q$-rep of lowest weight }\lambda-2\rho_l.
\]
\end{lemma}

\begin{proof}
Via duality and Lemma \ref{lem:1361} it suffices to show that the socle in $\Lind(\mu)$ is identified with the unique simple $\dot{\mfk{u}}_q$-representation of highest weight $\mu+ 2\rho_l$, at arbitrary $\mu$.  For this we note that $\Lind(\mu)=\mfk{v}_q^+$ as a $\mfk{v}^+_q$-representation.  So this representation has a unique highest weight vector of weight $\mu+2\rho_l$ if and only if $\mfk{v}_q^+$ has $1$-dimensional socle as a module over itself, i.e.\ $1$ dimensional space of $\mfk{v}_q^+$-invariants.
\par

Now, any invariant vector in $\mfk{v}_q^+$ decomposes into a sum of homogeneous invariant vectors, under the natural $Q$-grading on $\mfk{v}^+_q$, and any homogeneous invariant vector contributes a unique non-vanishing $\mfk{v}_q^{\geq 0}$-invariant vector in the free module
\[
\mfk{v}_q^{\geq 0}=\bigoplus_{\chi} \mfk{v}_q^{\geq 0}\ot_{k[\msf{K}_\alpha:\alpha\in \Delta]} k_\chi.
\]
Here the $\chi$ run over all characters for the (finite abelian) group generated by the units $\msf{K}_\alpha$ in $\mfk{v}_q$.  However, we know that the invariants in any finite-dimensional Hopf algebra are $1$-dimensional \cite[Theorem 2.1.3]{montgomery93} \cite{larsonsweedler69}.  So we have
\[
\dim \Lind(\mu)^{\mfk{v}^+_q}=\dim(\mfk{v}_q^+)^{\mfk{v}_q^+}\leq 1.
\]
\par

Since the unique vector $v$ of weight $\mu+2\rho_l$ in $\Lind(\mu)$ provides such an invariant vector, we conclude
\[
\Lind(\mu)^{\mfk{v}^+_q}=\Lind(\mu)_{\mu+2\rho}=k\cdot v.
\]
Uniqueness of the highest weight vector calculates the socle in $\Lind(\mu)$ as the simple $\dot{\mfk{u}}_q$-rep of highest weight $\mu+2\rho_l$.
\end{proof}

\begin{lemma}\label{lem:1396}
For any weight $\lambda\in P$, the identity map is the unique $\dot{\mfk{u}}_q$-linear endomorphism of $\Rind(\lambda)$, up to scaling.
\end{lemma}

\begin{proof}
Since $\Rind(\lambda)$ is of highest weight $\lambda$, and $\Rind(\lambda)_\lambda$ is $1$-dimensional, we have
\[
\End_{\dot{\mfk{u}}_q}(\Rind(\lambda))\cong\Hom_{\dot{\mfk{u}}_q^{\leq 0}}(\Rind(\lambda),k_\lambda)\cong\End_{\dot{\mfk{u}}_q^{\leq 0}}(k_\lambda)=k\cdot id.
\]
\end{proof}

\subsection{The Steinberg representation and $\dot{\mfk{u}}_q$-extensions}

We consider the Steinberg representation $\opn{St}$ in $\Rep(G_q)$, which we recall is the simple representation $\opn{St}=L(\rho_l)$ of highest weight $\rho_l=\frac{1}{2}\sum_{\gamma\in \Phi^+}(l_\gamma-1)\gamma$.  By Lemma \ref{lem:res_simples}, the restriction $\opn{St}|_{\dot{\mfk{u}}_q}$ to $\Rep(\dot{\mfk{u}}_q)$ is still simple.  It is, more specifically, the unique simple $\dot{\mfk{u}}_q$-representation of highest weight $\rho_l$.

\begin{lemma}[{cf.\ \cite[Theorem 9.10.3]{parshallwang91}}]\label{lem:1410}
At $\rho_l$ we have $\Rind(\rho_l)\cong\opn{St}|_{\dot{\mfk{u}}_q}$.  Furthermore, this object is self-dual.
\end{lemma}

\begin{proof}
We have $w_0\rho_l=-\rho_l$ so that $\opn{St}$ has lowest weight $-\rho_l$ and is self-dual, by Lemma \ref{lem:lowest_wt}.  Hence its restriction $\opn{St}|_{\dot{\mfk{u}}_q}$ is self-dual as well.  So we need only establish the claimed isomorphism $\Rind(\rho_l)\cong \opn{St}|_{\dot{\mfk{u}}_q}$.
\par

By Lemmas \ref{lem:1370} and \ref{lem:cosoc} there is a nontrivial inclusion $i:\opn{St}|_{\dot{\mfk{u}}_q}\to \Rind(\rho_l)$ and projection $p:\Rind(\rho_l)\to \opn{St}|_{\dot{\mfk{u}}_q}$.  Since these two representations have $1$-dimensional highest weight space of weight $\rho_l$, both $i$ and $p$ restrict to isomorphisms on the respective highest weight spaces.  In particular, the composites $ip$ and $pi$ are both nonzero.  Since
\[
\dim\End_{\dot{\mfk{u}}_q}(\opn{St}|_{\dot{\mfk{u}}_q})=\dim\End_{\dot{\mfk{u}}_q}(\Rind(\rho_l))=1
\]
we conclude that both of these composites are isomorphisms, and hence that both $i$ and $p$ are isomorphisms.
\end{proof}

\begin{lemma}
We have isomorphisms of $\dot{\mfk{u}}_q$-representations
\[
\Rind(\rho_l)\cong \Lind(-\rho_l)\cong \Rind_{\dot{\mfk{u}}^{\geq 0}_q}^{\dot{\mfk{u}}_q}(-\rho_l)\cong \Lind_{\dot{\mfk{u}}^{\geq 0}_q}^{\dot{\mfk{u}}_q}(\rho_l).
\]
\end{lemma}

\begin{proof}
By Lemmas \ref{lem:1361} and Lemma \ref{lem:1410} we have
\[
\Lind(-\rho_l)=\Rind(\rho_l)^\ast\cong \Rind(\rho_l).
\]
By exchanging $q$ with $q^{-1}$ and considering the algebra isomorphism $\psi:\dot{\mfk{u}}_q\overset{\sim}\to \dot{\mfk{u}}_{q^{-1}}$ which exchanges $E_\alpha$ with $F_\alpha$, and takes $1_\lambda$ to $1_{-\lambda}$, we see that
\[
\res_\psi\Rind_{\dot{\mfk{u}}_{q^{-1}}^{\leq 0}}^{\dot{\mfk{u}}_{q^{-1}}}(\rho_l)=\Rind_{\dot{\mfk{u}}_q^{\geq 0}}^{\dot{\mfk{u}}_q}(-\rho_l).
\]
Since $\opn{St}|_{\dot{\mfk{u}}_q}$ is simple with highest and lowest weights $\rho_l$ and $-\rho_l$ respectively we also have
\[
\opn{St}|_{\dot{\mfk{u}}_q}=\res_\psi(\opn{St}|_{\dot{\mfk{u}}_{q^{-1}}})=\res_\psi\Rind_{\dot{\mfk{u}}_{q^{-1}}^{\leq 0}}^{\dot{\mfk{u}}_{q^{-1}}}(\rho_l).
\]
So $\Rind(\rho_l)\cong \opn{St}|_{\dot{\mfk{u}}_q}\cong \Rind_{\dot{\mfk{u}}_q^{\geq 0}}^{\dot{\mfk{u}}_q}(-\rho_l)$.  Finally by fiddling with (self-)dualities again we obtain $\Rind_{\dot{\mfk{u}}_q^{\geq 0}}^{\dot{\mfk{u}}_q}(-\rho_l)\cong \Lind_{\dot{\mfk{u}}_q^{\geq 0}}^{\dot{\mfk{u}}_q}(\rho_l)$.
\end{proof}

We now observe vanishing of extensions for the Steinberg representation.

\begin{proposition}[{cf.\ \cite[Lemma 4.2]{andersenpolowen92}}]\label{prop:1466}
For any finite-dimensional $\dot{\mfk{u}}_q$-representation $V$ we have
\[
\Ext^1_{\dot{\mfk{u}}_q}(\opn{St}|_{\dot{\mfk{u}}_q},V)=\Ext^1_{\dot{\mfk{u}}_q}(V,\opn{St}|_{\dot{\mfk{u}}_q})=0.
\]
\end{proposition}

\begin{proof}
By duality it suffices to prove $\Ext^1_{\dot{\mfk{u}}_q}(V,\opn{St}|_{\dot{\mfk{u}}_q})=0$, and by considering exact sequences on cohomology it suffices to establish such vanishing at simple $V$.  Consider an extension
\begin{equation}\label{eq:1475}
0\to \opn{St}|_{\dot{\mfk{u}}_q}\to W\to V\to 0,
\end{equation}
with $V$ simple and of highest weight $\mu$.  Then, by Lemma \ref{lem:1370}, $V$ has lowest weight $\mu'$ with $\mu'\geq \mu-2\rho_l$.
\par

In the case where $\mu\ngeq \rho_l$, $W$ has a unique nonzero highest weight vector of weight $\rho_l$, up to scaling, and the $\dot{\mfk{u}}^{\leq 0}_q$-submodule $W'$ in $W$ generated by the all weight spaces $W_\tau$ with $\tau\neq \rho_l$ is of codimension $1$ with quotient $W/W'\overset{\sim}\to W_{\rho_l}\cong k_{\rho_l}$.  This projection lifts to a map $W\to \Rind(\rho_l)=\opn{St}|_{\dot{\mfk{u}}_q}$ which splits the sequence \eqref{eq:1475}.  A similar argument splits the sequence when $\mu=\rho_l$.
\par

In the final case where $\mu >\rho_l$, $W$ has $1$-dimensional lowest weight space of weight $-\rho_l$.  We have the $\dot{\mfk{u}}^{\geq 0}_q$-projection $W\to W_{-\rho_l}\cong k_{-\rho_l}$ and corresponding map $W\to \Rind_{\dot{\mfk{u}}^{\geq 0}_q}^{\dot{\mfk{u}}_q}(-\rho_l)=\opn{St}|_{\dot{\mfk{u}}_q}$, which again splits \eqref{eq:1475}.
\end{proof}

\subsection{Projective and injective $\dot{\mfk{u}}_q$-representations}

\begin{lemma}\label{lem:1490}
The category $\rep(\dot{\mfk{u}}_q)$ of finite-dimensional $\dot{\mfk{u}}_q$-representations has enough injectives, and any injective object in $\rep(\dot{\mfk{u}}_q)$ is injective in $\Rep(\dot{\mfk{u}}_q)$.
\end{lemma}

\begin{proof}
As in \cite[pg 118]{parshallwang91} we consider the right adjoint $I:\opn{Vect}_P=\Rep(T_q)\to \Rep(\dot{\mfk{u}}_q)$ to the restriction functor $\Rep(\dot{\mfk{u}}_q)\to \Rep(T_q)$.  We note that $I$ sends injectives to injectives since restriction is exact.  By an explicit expression
\[
I(V)=\varinjlim_\alpha\Hom_{T_q}(\dot{\mfk{u}}_q,V_\alpha)\cong \varinjlim_\alpha\Hom_k(\mfk{v}_q^+\ot \mfk{v}_q^-,V_\alpha)\cong \Hom_k(\mfk{v}_q^+\ot \mfk{v}_q^-,V),
\]
where the colimit here is indexed over finite-dimensional subrepresentations $V_\alpha\subseteq V$, we see that $I$ preserves finite-dimensional objects.
\par

Now, since all objects in $\Rep(T_q)$ are injective, the unit of the adjunction $V\to I(V)$ maps any finite-dimensional $\dot{\mfk{u}}_q$-representation to a finite-dimensional injective object.  Since the restriction functor is faithful, this unit map $V\to I(V)$ has trivial kernel.  So we see that $\rep(\dot{\mfk{u}}_q)$ has enough injectives, and injectives in $\rep(\dot{\mfk{u}}_q)$ are injective in $\Rep(\dot{\mfk{u}}_q)$.
\end{proof}

\begin{proposition}\label{prop:st_dotu}
The Steinberg representation $\opn{St}$ restricts to a simultaneous projective and injective object in $\Rep(\dot{\mfk{u}}_q)$.
\end{proposition}

\begin{proof}
Proposition \ref{prop:1466} tells us that $\opn{St}|_{\dot{\mfk{u}}_q}$ is injective in $\rep(\dot{\mfk{u}}_q)$.  By Lemma \ref{lem:1490}, it follows that $\opn{St}|_{\dot{\mfk{u}}_q}$ is injective in $\Rep(\dot{\mfk{u}}_q)$.  For projectivity, Proposition \ref{prop:1466} says that $\opn{St}|_{\dot{\mfk{u}}_q}$ is also projective in $\rep(\dot{\mfk{u}}_q)$.  For a possibly infinite-dimensional $\dot{\mfk{u}}_q$-representation $M$, we can write $M$ as a filtered colimit $M=\varinjlim_\alpha M_\alpha$ of finite-dimensional subrepresentations.  We then employ compactness of $\opn{St}|_{\dot{\mfk{u}}_q}$ in $\Rep(\dot{\mfk{u}}_q)$ and exactness of filtered colimits to find
\[
\Ext^1_{\dot{\mfk{u}}_q}(\opn{St}|_{\dot{\mfk{u}}_q},M)=\varinjlim_\alpha \Ext^1_{\dot{\mfk{u}}_q}(\opn{St}|_{\dot{\mfk{u}}_q},M_\alpha)\underset{\rm Prop\ \ref{prop:1466}}=0.
\]
Such vanishing ensures projectivity of $\opn{St}|_{\dot{\mfk{u}}_q}$ in $\Rep(\dot{\mfk{u}}_q)$.
\end{proof}

\begin{corollary}\label{cor:1627}
The categories $\rep(\dot{\mfk{u}}_q)$ and $\Rep(\dot{\mfk{u}}_q)$ have enough projectives, and projectives in $\rep(\dot{\mfk{u}}_q)$ remain projective in $\Rep(\dot{\mfk{u}}_q)$.
\end{corollary}

\begin{proof}
Take $P=(\opn{St}|_{\dot{\mfk{u}}_q})\ot(\opn{St}|_{\dot{\mfk{u}}_q})=(\opn{St}|_{\dot{\mfk{u}}_q})^\ast\ot(\opn{St}|_{\dot{\mfk{u}}_q})$.  This module is projective and coevaluation provides a surjection $P\to \1$.  Tensoring with any representation $V$ provides a surjection $P\ot V\to V$, and the object $P\ot V$ is seen to be projective via the existence of the adjunction
\[
\Hom_{\dot{\mfk{u}}_q}(P\ot V,-)\cong \Hom_{\dot{\mfk{u}}_q}(V,{^\ast P}\ot-)\cong \Hom_{\dot{\mfk{u}}_q}(V,P\ot -).
\]
\par

To elaborate, $P\ot -$ is an exact functor which takes injective values in $\Rep(\dot{\mfk{u}}_q)$ (see Lemma \ref{lem:inj_A}).  Hence $P\ot-$ sends exact sequences to split exact sequences, and we observe exactness of the above functor.  We note finally that $P\ot V$ is finite-dimensional whenever $V$ is finite-dimensional to see that $\rep(\dot{\mfk{u}}_q)$ has enough projectives.
\end{proof}

\section{Restriction, projectivity, and injectivity for simply-connected $G$}
\label{sect:steinborg_sc}

We prove that the Steinberg representation is both projective and injective in $\Rep(G_q)$, in the simply-connected setting, and that it restricts to a simultaneously projective and injective representation in $\Rep(\bar{u}_q)$.  We use this fact to find that the restriction map
\[
\res^-_q:\Rep(G_q)\to \Rep(\bar{u}_q)
\]
both preserves, and detects, projective and injective objects in $\Rep(G_q)$.

\subsection{Vibe check}

In Section \ref{sect:steinborg} we considered two ``induction" functors, the left and right adjoints to restriction $\Rep(\dot{\mfk{u}}_q)\to \Rep(\dot{\mfk{u}}^{\leq 0}_q)$ for torally extended $\mfk{v}_q$-representations.  We adopted the notations $\Lind$ and $\Rind$ for these functors.
\par

From this point on we only consider \emph{right adjoints} to restriction functors, and any functor which we refer to as \emph{induction} is the right adjoint to some restriction functor.

\subsection{Exactness of induction}

We first recall a standard fact.

\begin{lemma}[{\cite[\S 3.3]{etingofostrik04}}]\label{lem:1699}
Let $A$ be a Hopf algebra.  For any map of $A$-module coalgebras $A\to C$, the right adjoint to restriction $\opn{ind}:\opn{Corep}(C)\to \opn{Corep}(A)$ is a map of $\opn{Corep}(A)$-module categories.
\end{lemma}

\begin{proof}
Since the trivial corepresentation $\1$ is compact in $\opn{Corep}(C)$, the functor
\[
\opn{ind}(-)=\Hom_C(\1,-\ot A)=(-\ot A)^C=-\square_C A
\]
\cite[Proposition 6]{doi81} commutes with filtered colimits.  So it suffices to establish the proposed natural isomorphism $\opn{ind}(-\ot\res(V))\cong \opn{ind}(-)\ot V$ at finite-dimensional $V$.  But now the claim follows by duality, since we have the sequence of natural identifications
\[
\begin{array}{rl}
\Hom_C(\res(-),-\ot \res(V))& \cong\Hom_C(\res(-\ot {^\ast V}),-)\\
&\cong\Hom_A(-\ot {^\ast V},\opn{ind}(-))\cong\Hom_A(-,\opn{ind}(-)\ot V).
\end{array}
\]
\end{proof}

Via Tannakian reconstruction, Lemma \ref{lem:1699} tells us that the right adjoints to the restriction functors
\[
\Rep(G_q)\to \Rep(\dot{\mfk{u}}_q),\ \ \Rep(\dot{\mfk{u}}_q)\to \Rep(\bar{u}_q),\ \text{and}\ \ \Rep(G_q)\to \Rep(\bar{u}_q)
\]
are naturally $\Rep(G_q)$-linear, $\Rep(\dot{\mfk{u}}_q)$-linear, and again $\Rep(G_q)$-linear, respectively.

\begin{proposition}\label{prop:1234}
Let $G$ be simply-connected.  The right adjoint to restriction $\Rep(\bar{u}_q)\to \Rep(\dot{\mfk{u}}_q)$ is faithfully exact.
\end{proposition}

\begin{proof}
We can just write down the adjoint.  For any $\bar{u}_q$-representation $V$ there is a unique lift of $V$ to a $P$-graded vector space $\opn{lift}(V)$ whose $P$-grading is supported on $P_l$.  We then take the orbit under the action of the simples $k_\lambda$ associated to $\lambda\in P^\ast$ to get the $P$-graded vector space
\[
\opn{ind}(V)=\bigoplus_{\lambda\in P^\ast}k_\lambda\ot \opn{lift}(V).
\]
This graded vector space carries a natural $\dot{\mfk{u}}_q$-action specified by the formulae
\[
\mbf{E}_\alpha\cdot (\zeta_\lambda\ot v)=\zeta_{\lambda'}\ot (\mbf{E}_\alpha\cdot v)\ \ \text{and}\ \ F_\alpha\cdot (\zeta_\lambda\ot v)=\zeta_{\lambda''}\ot (F_\alpha\cdot v),
\]
where 
\[
\lambda'=\lambda+\opn{deg}(v)+\alpha-\deg(\mbf{E}_\alpha\cdot v)\ \ \text{and}\ \ \lambda''=\lambda+\opn{deg}(v)-\alpha-\deg(F_\alpha\cdot v).
\]
This operation is functorial, is seen directly to provide a right adjoint to restriction, and is faithful and exact by construction.
\end{proof}

\begin{proposition}\label{prop:1247}
Let $G$ be simply-connected.  The right adjoint to restriction $\Rep(\dot{\mfk{u}}_q)\to \Rep(G_q)$ is exact.
\end{proposition}

We proceed as in the proof of \cite[Theorem 4.8]{andersenpolowen91}.

\begin{proof}
Let $\1\overset{\sim}\to E$ be a bounded below resolution of $\1$ in $\Rep(G_q)$ by objects of the form $E^i=\bar{E}^i\ot \opn{St}$, where $\bar{E}^i$ is any representation in $\Rep(G_q)$.  Then $E$ restricts to an injective resolution of the unit $\1\to E|_{\dot{\mfk{u}}_q}$ in $\Rep(\dot{\mfk{u}}_q)$.  For any $V$ in $\Rep(\dot{\mfk{u}}_q)$ we have the injective resolution $V\ot E|_{\dot{\mfk{u}}_q}$ and calculate the higher derived functors of induction as
\[
\begin{array}{rl}
\opn{ind}^{>0}(V) & = H^{>0}\big(\opn{ind}(V\ot E|_{\dot{\mfk{u}}_q})\big)\\
& \cong H^{>0}\big(\opn{ind}(V)\ot E\big)=\opn{ind}(V)\ot H^{>0}(E)=0.
\end{array}
\]
Since the higher derived functors vanish, $\opn{ind}$ is exact.
\end{proof}

We combine Propositions \ref{prop:1234} and \ref{prop:1247} to see that induction from the smallest quantum algebra to $\Rep(G_q)$ is also exact.

\begin{theorem}\label{thm:ind_ex}
Let $G$ be simply-connected.  The right adjoint to restriction $\Rep(\bar{u}_q)\to \Rep(G_q)$ is faithfully exact.
\end{theorem}

\begin{proof}
The restriction functor $\res^-_q$ factors through $\Rep(\dot{\mfk{u}}_q)$.  Hence the right adjoint to $\res^-_q$ factors as
\[
\Rep(\bar{u}_q)\overset{\opn{ind}_1}\longrightarrow \Rep(\dot{\mfk{u}}_q)\overset{\opn{ind}_2}\longrightarrow \Rep(G_q).
\]
So exactness follows immediately from Propositions \ref{prop:1234} and \ref{prop:1247}.
\par

Via exactness, we can prove faithfulness by observing non-vanishing of $\opn{ind}(L)$ for each simple $L$ in $\Rep(\bar{u}_q)$.  However, each simple $L$ in $\Rep(\bar{u}_q)$ is the image of some simple $L'$ in $\Rep(G_q)$, by Corollary \ref{cor:simples}.  So we have
\[
0\neq \Hom_{\bar{u}_q}\big(\res^-_q(L'),L\big)=\Hom_{\bar{u}_q}\big(L',\opn{ind}(L)\big).
\]
In particular, $\opn{ind}(L)$ is nonzero.
\end{proof}

\begin{remark}
For $q$ with odd order scalar parameters, Theorem \ref{thm:ind_ex} appears in \cite[Theorem 4.8]{andersenpolowen92}.
\end{remark}

\subsection{Projectives and injectives}

\begin{proposition}\label{prop:injectives}
Let $G$ be simply-connected.  A $G_q$-representation is injective in $\Rep(G_q)$ if and only if its image in $\Rep(\bar{u}_q)$ is injective.
\end{proposition}

\begin{proof}
First note that the induction functor $\opn{ind}:\Rep(\bar{u}_q)\to \Rep(G_q)$ satisfies
\[
\opn{ind}(\1)=\1\square_{(\bar{u}_q)^\ast}\O(G_q)={^{\bar{u}_q^\ast}\O(G_q)}.
\]
By Corollary \ref{cor:coinv} this object is just $\O(G^\ast_{\varepsilon})$, and in particular lies in the semisimple subcategory $\Rep(G^\ast_\varepsilon)=\opn{ker}(\res^-_q)$ in $\Rep(G_q)$.  Hence the unit map of the adjunction $\1\to \opn{ind}(\1)$ is split.  It follows that $\opn{ind}(\1)$ contains $\1$ as a summand, and that each $G_q$-representation $V$ appears as a summand in the induction $\opn{ind}(\res^-_q V)\cong \opn{ind}(\1)\ot V$.
\par

Suppose first that $V$ has injective image in $\Rep(\bar{u}_q)$.  Since $V$ appears as a summand in $\opn{ind}(\1)\ot V$, it suffices to show that $\opn{ind}(\1)\ot V$ is injective in $\Rep(G_q)$.  We have the natural identification
\[
\Hom_{G_q}(-,\opn{ind}(\1)\ot V)\cong\Hom_{G_q}(-,\opn{ind}(\res^-_qV))\cong\Hom_{\bar{u}_q}(\res^-_q-,\res^-_qV),
\]
from which we conclude exactness of the functor $\Hom_{G_q}(-,\opn{ind}(\1)\ot V)$.  It follows that $\opn{ind}(\1)\ot V$ is injective, and hence that $V$ is injective.
\par

Suppose conversely that $V$ is injective in $\Rep(G_q)$.  To verify injectivity of $\res^-_q(V)$ in $\Rep(\bar{u}_q)$, we show that the extensions $\Ext^{>0}_{\bar{u}_q}(L,\res^-_q V)$ vanish for each simple $L$ in $\Rep(\bar{u}_q)$.  By Corollary \ref{cor:simples}, we can find a representation $L'$ in $\Rep(G_q)$ with $\res^-_q(L')=L$, and by exactness of induction we obtain an identification
\begin{equation}\label{eq:1755}
\Ext^{>0}_{\bar{u}_q}(L,\res^-_qV)\cong \Ext^{>0}_{G_q}(L',\opn{ind}(\res^-_q V))\cong \Ext^{>0}_{G_q}(L',\opn{ind}(\1)\ot V).
\end{equation}
Since $V$ is injective, the product $\opn{ind}(\1)\ot V$ is injective in $\Rep(G_q)$.  Hence the extensions \eqref{eq:1755} vanish.
\end{proof}

For projectivity, we see most clearly that one can detect \emph{finite-dimensional} projectives via restriction, though we remove this constraint in a moment.

\begin{lemma}\label{lem:fd_projectives}
The restriction functor $\Rep(G_q)\to \Rep(\bar{u}_q)$ preserves, and detects, finite-dimensional projective objects.
\end{lemma}

\begin{proof}
Suppose that $V$ is projective in $\Rep(G_q)$.  Then the functor
\[
\Hom_{\bar{u}_q}(\res^-_q(V),-)=\Hom_{G_q}(V,\opn{ind}(-))
\]
is exact by Theorem \ref{thm:ind_ex}.  So $\res^-_q(V)$ is projective in $\Rep(\bar{u}_q)$.

Suppose conversely that $V$ is finite-dimensional and has projective image in $\Rep(\bar{u}_q)$.  We have an identification
\[
\Hom_{G_q}(V,-)=\Hom_{G^\ast_{\varepsilon}}(\1,\Hom_{\bar{u}_q}(\res^-_q V,-)).
\]
(Note that we employ finiteness of $V$ here to ensure that, at any $G_q$-representation $W$, the $G^\ast_\varepsilon$-representation $\Hom_{\bar{u}_q}(\res^-_qV,W)$ is integrable.)  This functor is exact by projectivity of $\res^-_q(V)$ and semisimplicity of the category $\Rep(G^\ast_\varepsilon)$.
\end{proof}

By completely similar arguments, with $G^\ast_\varepsilon$ and Proposition \ref{prop:norm} replaced by $T^\ast_\varepsilon$ and Lemma \ref{lem:norm_dot}, respectively, we obtain toral analogs of Proposition \ref{prop:injectives} and Lemma \ref{lem:fd_projectives}.

\begin{proposition}\label{prop:1801}
An object $V$ in $\Rep(\dot{\mfk{u}}_q)$ is injective if and only if its image in $\Rep(\bar{u}_q)$ is injective.  A finite-dimensional object $V$ in $\Rep(\dot{\mfk{u}}_q)$ is projective if and only if its image in $\Rep(\bar{u}_q)$ is projective.
\end{proposition}

We recall that the Steinberg representation is projective and injective in $\Rep(\dot{\mfk{u}}_q)$, by Proposition \ref{prop:st_dotu}.  So Proposition \ref{prop:1801} implies projectivity and injectivity of $\opn{St}$ over $\bar{u}_q$.  We recall also that $\opn{St}$ has simple restriction to $\bar{u}_q$, by Proposition \ref{prop:res_simples}.

\begin{corollary}\label{cor:1757}
For simply-connected $G$, the Steinberg representation $\opn{St}$ is projective, injective, and simple over $\bar{u}_q$.
\end{corollary}

We now ``go backwards" to observe projectivity and injectivity of $\opn{St}$ over $G_q$.

\begin{corollary}[\cite{andersen03}]
Let $G$ be simply-connected.  The Steinberg representation $\opn{St}$ is both projective and injective in $\Rep(G_q)$, the subcategory $\rep(G_q)$ of finite-dimensional $G_q$-representations has enough projectives, and any projective in $\rep(G_q)$ is projective in $\Rep(G_q)$.
\end{corollary}

\begin{proof}
The first claim follows by the three previous results.  The subsequent claims are established just as in the proof of Corollary \ref{cor:1627}.
\end{proof}

We now understand infinite-dimensional projectives in $\Rep(G_q)$ via the following observation.

\begin{lemma}\label{lem:1395}
Let $G$ be simply-connected.  An object $V$ in $\Rep(G_q)$ is projective if and only if $V$ is a summand of $V\ot \opn{St}\ot\opn{St}$.  Similarly, $V$ is injective if and only if $V$ is a summand of $V\ot \opn{St}\ot \opn{St}$.
\end{lemma}

\begin{proof}
If $I$ is an injective representation then $W\ot I$ is injective for all $W$ (see Lemma \ref{lem:inj_A}).  So all products $W\ot \opn{St}$ are injective.  Furthermore, tensoring with the coevaluation map $\1\to \opn{St}\ot \opn{St}^\ast=\opn{St}\ot \opn{St}$ provides an injection $V\to V\ot \opn{St}\ot \opn{St}$ at arbitrary $V$.  In the case that $V$ is injective, this injection splits, so that $V$ is a summand of the given product.  Furthermore, since the class of injectives is closed under taking summands, we see that this property characterizes injectives.
\par

For projectives, we first claim that $W\ot \opn{St}$ is projective for any $W$.  This is clear from self-duality of $\opn{St}$, the adjunction
\[
\Hom_{G_q}(W\ot \opn{St},-)\cong \Hom_{G_q}(W,-\ot \opn{St}^\ast)\cong \Hom_{G_q}(W,-\ot\opn{St}),
\]
and the fact that $-\ot\opn{St}$ is an exact endomorphism to the subcategory of injectives in $\Rep(G_q)$.  Given this fact, we have the evaluation $\opn{St}\ot \opn{St}=\opn{St}^\ast\ot \opn{St}\to \1$ which provides a surjection $V\ot \opn{St}\ot \opn{St}\to V$ at arbitrary $V$.  When $V$ is projective, this surjection splits, and conversely splitting of this map implies projectivity of $V$.  So we establish the claimed result.
\end{proof}

As a corollary we see that the category $\Rep(G_q)$ is Frobenius, in the sense that its projectives and injectives agree.

\subsection{Simply-connected conclusions}

We collect our findings.

\begin{theorem}\label{thm:steinberg}
Let $G$ be simply-connected.  Consider the dominant weight $\rho_l=\frac{1}{2}\sum_{\gamma\in \Phi^+}(l_\gamma-1)\gamma$ and the corresponding simple $G_q$-representation $\opn{St}=L(\rho_l)$.\vspace{1mm}
\begin{enumerate}
\item The representation $\opn{St}$ is both projective and injective in $\Rep(G_q)$.\vspace{2mm}
\item The representation $\opn{St}$ restricts to an object which is simultaneously projective, injective, and simple in $\Rep(\bar{u}_q)$.\vspace{2mm}
\item The category $\Rep(G_q)$ has enough projectives and injectives, and is Frobenius.\vspace{2mm}
\item The subcategory $\rep(G_q)$ of finite-dimensional $G_q$-representations has enough projectives and injectives, and is Frobenius.\vspace{2mm}
\item The inclusion $\rep(G_q)\to \Rep(G_q)$ preserves (and detects) projective and injective objects.\vspace{2mm}
\item An object in $\Rep(G_q)$ is projective and injective if and only if its restriction to $\Rep(\bar{u}_q)$ is projective and injective.
\end{enumerate}
\end{theorem}

\section{Restriction, projectivity, and injectivity for general $G$}
\label{sect:steinborg_nsc}

In this final section we generalize Theorem \ref{thm:steinberg} to arbitrary semisimple groups.

\subsection{Homological algebra along the quotient $G^{sc}_q\to G_q$}
\label{sect:2213}

\begin{lemma}\label{lem:2215}
Let $G$ be a semisimple algebraic group, and let $G^{sc}$ be the simply-connected form for $G$.
\begin{enumerate}
\item The categories $\rep(G_q)$ and $\Rep(G_q)$ have enough projectives and injectives.\vspace{2mm}
\item The inclusion functor $\rep(G_q)\to \Rep(G_q)$ preserves (and detects) projective and injective objects.\vspace{2mm}
\item An object in $\Rep(G_q)$ is projective if and only if it is injective.\vspace{2mm}
\item An object in $\Rep(G_q)$ is projective and injective if and only if its image in $\Rep(G^{sc}_q)$ is projective and injective.
\end{enumerate}
\end{lemma}

\begin{proof}
(4) Note that $\Rep(G_q)$ is equal to the full subcategory of representations in $\Rep(G^{sc}_q)$ whose $P$-grading is supported on the character lattice $X$ for $G$.  We let $F:\Rep(G_q)\to \Rep(G^{sc}_q)$ denote the corresponding monoidal embedding and consider the group of characters $\Sigma=(P/X)^\vee$.  This group acts naturally on $G^{sc}_q$-representations, and the invariants $W^{\Sigma}$ in any $G^{sc}_q$-representation $W$ produce a maximal $G_q$-subrepresentation in $W$.  Indeed, taking $\Sigma$-invariants simply isolates the subrepresentation in $W$ consisting of those vectors which are supported on the sublattice $X$ in $P$.
\par

We have the natural $G^{sc}_q$-linear inclusion $i_W:W^\Sigma\to W$ and projection $p_W:W\to W^\Sigma$, and restrictions along $i_W$ and $p_W$ provide adjunctions
\begin{equation}\label{eq:1905}
\Hom_{G^{sc}_q}(F-,W)=\Hom_{G_q}(-,W^{\Sigma})\ \ \text{and}\ \ \Hom_{G^{sc}_q}(W,F-)=\Hom_{G_q}(W^{\Sigma},-).
\end{equation}
Since taking $\Sigma$-invariants is an exact functor, the above adjunctions imply that a $G_q$-representation $P$ has projective (resp.\ injective) image in $\Rep(G^{sc}_q)$ whenever $P$ is projective (resp.\ injective) in $\Rep(G_q)$.

Conversely, fully faithfulness of the embedding $F$ implies that a given representation $P$ over $G_q$ is projective (resp.\ injective) in $\Rep(G_q)$ provided it has projective (resp.\ injective) image in $\Rep(G^{sc}_q)$.  We therefore establish (4).
\par

(1) Consider any $G_q$-representation $V$ and any surjection $P\to FV$ from a projective in $\Rep(G^{sc}_q)$.  By Theorem \ref{thm:steinberg}, such $P$ exists and can be taken to be finite-dimensional when $V$ is finite-dimensional.  By the adjunction \eqref{eq:1905} we find that the invariants $P^{\Sigma}$ are projective in $\Rep(G_q)$, and the induced surjection $P^{\Sigma}\to V$ provides a projective covering of $V$.  One similarly argues with injectives, so that we obtain (1).  Properties (2) and (3) are now inherited from $\Rep(G^{sc}_q)$ via (1), (4), and Theorem \ref{thm:steinberg}.
\end{proof}

For the smallest quantum algebra, we have the inclusion $X\to P$ and recover $X^\ast$ as the intersection $X^\ast=X\cap P^\ast$.  So this inclusion reduces to an inclusion on the quotients $X/X^\ast\to P/P^\ast$.  We therefore locate $\Rep(\bar{u}_q)$ in its simply-connected form $\Rep(\bar{u}_q^{sc})$ as the full subcategory of $\bar{u}_q^{sc}$-representations whose $P/P^\ast$-grading is supported on the subgroup $X/X^\ast$.

\begin{lemma}\label{lem:1456}
Consider the smallest quantum algebra $\bar{u}_q^{sc}$ for $G^{sc}_q$, and the inclusion $\Rep(\bar{u}_q)\to \Rep(\bar{u}^{sc}_q)$.  A $\bar{u}_q$-representation is projective (resp.\ injective) if and only if its image in $\Rep(\bar{u}^{sc}_q)$ is projective (resp.\ injective).
\end{lemma}

\begin{proof}
Since the inclusion $\Rep(\bar{u}_q)\to \Rep(\bar{u}_q^{sc})$ is exact and fully faithful, any $\bar{u}_q$-representation $V$ whose image in $\Rep(\bar{u}_q^{sc})$ is projective (resp.\ injective) must be projective (resp.\ injective) in $\Rep(\bar{u}_q)$.
\par

For the converse, take $\opn{K}^\vee=\opn{coker}(X/X^\ast\to P/P^\ast)$ and consider the central subgroup $\opn{K}=(\opn{K}^{\vee})^{\vee}\subseteq \bar{u}^{sc}_q$.  For any $\bar{u}_q$-representation $V$ and $\bar{u}^{sc}_q$-representation $W$, we have natural identifications
\[
\Hom_{\bar{u}^{sc}_q}(V,W)=\Hom_{\bar{u}_q}(V,W^{\opn{K}})\ \ \text{and}\ \ \Hom_{\bar{u}^{sc}_q}(W,V)=\Hom_{\bar{u}_q}(W^{\opn{K}},V).
\]
These identifications are obtained by composing with the natural inclusion $W^{\opn{K}}\to W$ and projection $W\to W^{\opn{K}}$ respectively.  Since $\opn{K}$ is finite, taking $\opn{K}$-invariants is an exact functor, and we conclude that $V$ is projective (resp.\ injective) in $\Rep(\bar{u}^{sc}_q)$ whenever $V$ is projective (resp.\ injective) in $\Rep(\bar{u}_q)$.
\end{proof}

\subsection{Semisimple algebraic groups}

We apply Lemmas \ref{lem:2215} and \ref{lem:1456} to obtain a version of Theorem \ref{thm:steinberg} for general $G$.

\begin{theorem}\label{thm:steinberg_semis}
Let $G$ be a semisimple algebraic group and $q$ be an arbitrary quantum parameter.  The following hold:\vspace{1mm}
\begin{enumerate}
\item The category $\Rep(\bar{u}_q)$ is Frobenius.\vspace{2mm}
\item The category $\Rep(G_q)$ has enough projectives and injectives, and is Frobenius.\vspace{2mm}
\item The subcategory $\rep(G_q)$ of finite-dimensional representations has enough projectives and injectives, and is Frobenius.\vspace{2mm}
\item The inclusion $\rep(G_q)\to \Rep(G_q)$ preserves (and detects) projective and injective objects.\vspace{2mm}
\item An object in $\Rep(G_q)$ is projective and injective if and only if its restriction to $\Rep(\bar{u}_q)$ is projective and injective.\vspace{2mm}
\end{enumerate}
\end{theorem}

\begin{proof}
Points (2)--(4) were already covered in Lemma \ref{lem:2215}.   Point (5) follows from Theorem \ref{thm:steinberg} (6) in conjunction with Lemma \ref{lem:2215} (5) and Lemma \ref{lem:1456}.  Point (1) follows from (2), (5), and Corollary \ref{cor:simples} in the simply-connected case.  The general case then follows by Lemma \ref{lem:1456}.  Alternatively, we observe (1) via Skryabin's Theorem \cite[Theorem 6.1]{skryabin07}.
\end{proof}

\subsection{An incomplete accounting of the literature}
\label{sect:history_inc}

As we noted in the preamble to Section \ref{sect:steinborg}, the Steinberg representation has a long and storied history in both quantum and modular representation theory.  We recall some of the quantum hits below.
\par

In \cite{parshallwang91} Parshall and Wang show that, for $G=\opn{SL}(n)$ and $q$ of odd order, the Steinberg representation is projective and injective in $\Rep(\SL(n)_q)$, and also in $\Rep(\bar{u}_q)$ \cite[Theorem 9.10.3, Corollary 9.10.4]{parshallwang91}.  Their proof relies on an analysis of a gaggle of induction functors, between various quantum groups, and a reciprocity theorem which essentially counts the lengths of indecomposable injectives in $\Rep(\dot{\mfk{u}}_q)$ \cite[Proposition 9.8.4]{parshallwang91}.  As the authors note, their methods are directly inspired by known phenomena in modular representation theory.  See for example \cite[Ch 10 \& 11]{jantzen03}.
\par

In \cite{andersenpolowen91} Andersen, Polo, and Wen prove that the Steinberg representation is projective and injective in $\Rep(G_q)$, but now for arbitrary simply-connected $G$ and $q$ of a prime power order \cite[Lemma 6.6, Theorem 9.8]{andersenpolowen91}.  In this context the authors rely on a rigorous analysis of induction from the quantum Borel $H^\ast(G_q/B_q,-)$, and a linkage principle for $G_q$-representations \cite[\S\ 8]{andersenpolowen91}.  In a subsequent text \cite{andersenpolowen92} Andersen, Polo, and Wen prove projectivity and injectivity of $\opn{St}$ over the small quantum group, under essentially the same restrictions on $G$ and $q$ \cite[Theorem 4.3]{andersenpolowen92}.  In \cite{andersenpolowen92} it is also shown that the induction functor from the small quantum group $\Rep(\bar{u}_q)\to \Rep(G_q)$ is exact and faithful \cite[Theorem 4.8]{andersenpolowen92}.

In \cite{andersenwen92} Andersen and Wen again work at $q$ of prime power order and provide analyses which parallel those of \cite{andersenpolowen91,andersenpolowen92}, but now in mixed characteristic.  In \cite{andersen92} Andersen addresses induction and tilting modules for simply-connected $G$ at $q$ of odd order which is greater than the Coxeter number.  In \cite{andersen03} Andersen returns to this topic, and proves a linkage principle for the big quantum group $G_q$ at simply-connected $G$ and arbitrary $q$.  As suggested in \cite[pg 14]{andersen03}, one can argue from the strong linkage principle that $\opn{St}$ is both projective and injective in $\Rep(G_q)$ \cite{andersen}, though the details are not entirely trivial and do not appear in the referenced text.
\par

The papers \cite{andersenpolowen92,andersenwen92,andersen03} provide the most recent treatments of these topics which are relatively general and completely rigorous, as far as I know.  Since \cite{andersenpolowen92,andersenwen92} it has often been proposed, either implicitly or explicitly, that the methods of these texts generalize immediately to most quantum parameters, most odd quantum parameters, or all parameters $q$.  (Consider, for example, the curious appearance of the Steinberg representation in \cite[proof of Lemma 4.1]{negron21}.)
\par

In reconsidering this topic however, the transition from the odd order case to the even order case did not seem to be completely transparent.  Hence we produce this text.

\appendix

\section{Braid group action on $\mfk{v}_q$}
\label{sect:B_g}

\emph{For this section $\mfk{v}_q$ specifically denotes the subalgebra in $U_q$ generated by the simple root vectors $E_\alpha$ and $F_\alpha$ with $l_\alpha>1$, the grouplikes $\msf{K}_{\nu}$, and the distinguished generators from Propositions \ref{prop:863} and \ref{prop:865} at low order $q$.}

Our goal is to show that the subalgebra $\mfk{v}_q$ is preserved under the action of the braid group operators on $U_q$.  It will follow that $\mfk{v}_q$ is the smallest subalgebra in $U_q$ which contains all of the simple root vectors $E_\alpha$ and $F_\alpha$ with $l_\alpha>1$, all of the grouplikes $\msf{K}_\nu$, and is stable under the action of the braid group $B_{\mfk{g}}$, i.e.\ that $\mfk{v}_q$ is identified with the subalgebra of \eqref{eq:mfkv_q}.
\par

Let us be clear here that our braid group $B_{\mfk{g}}$, for $\mfk{g}=\opn{Lie}(G)$, acts on $U_q$ via the specific operators
\[
T_\alpha=T''_{\alpha,1}
\]
from \cite[\S 37]{lusztig90II}.

\begin{remark}
The reader should note that the usual braid group \cite{artin47} (only) appears in type $A_n$, from the perspective of the works \cite{lusztig90II,lusztig93}.  The ``braid group" we are referring to here is a group $B_\mfk{g}$ which depends on the Dynkin type of $\mfk{g}$.  Specifically, the operators $T_\alpha$ in $B_{\mfk{g}}$ are only supposed to satisfy the relations
\[
T_{\alpha_m}\dots T_{\alpha_1}=T_{\beta_m}\dots T_{\beta_1}\ \ \text{whenever $\prod_i\sigma_{\alpha_i}=\prod_i\sigma_{\beta_i}$ and }\opn{length}(\sigma)=m.
\]
\end{remark}

\subsection{Reducing to almost-simple $G$ at the lacing number}

When each $l_i$ is greater than the lacing number for its associated almost-simple factor in $G$, the algebra $\mfk{v}_q$ is seen to be stable under the action of $B_{\mfk{g}}$ by direct inspection.  Namely, for each pair of simple roots $\alpha$ and $\beta$ the explicit expression of $T_\beta(E_\alpha)$ from \cite{lusztig90II} shows that $T_\beta(E_\alpha)$ is in the subalgebra generated by the simple $E_\nu$ in $U_q$.  The same is true of the $F$'s.
\par

According to this information, we need only check the case where some of the $l_i$'s are less than or equal to the lacing numbers for their associated factors in $G$.  The problem then reduces to considering almost-simple $G$ at such low order parameters.  We deal with these particular almost-simple cases below.

\subsection{Some lemmas}

\begin{lemma}\label{lem:T_inv}
If $\beta$ is a simple root with $q_\beta=-1$, then $T_\beta=T_\beta^{-1}$.
\end{lemma}

\begin{proof}
In this case $\msf{K}_\beta=\msf{K}_\beta^{-1}$, and $\msf{K}_\beta$ centralizes both $E_\beta$ and $F_\beta$.  So the result follows from the explicit expression of $T^{-1}_\beta$ \cite[\S 37.1.3]{lusztig93} and the fact that $(-1)^rq_\beta^{\pm r}=1$ at every integer index $r$.
\end{proof}

We translate between computations for the $E$'s and $F$'s via a ``symmetry" on the quantum enveloping algebra.  Specifically, we observe the following braid group equivariant isomorphism.

\begin{lemma}[\cite{lusztig90II}]
At an arbitrary parameter $q$, there is an anti-algebra isomorphism $\phi:U_{q^{-1}}\to U_{q}$,
\[
\phi(E_\alpha^{(n)})=F^{(n)}_\alpha,\ \phi(F^{(n)}_\alpha)=E_\alpha^{(n)},\ \phi(\msf{K}_\gamma)=\msf{K}_{\gamma}^{-1},
\]
with $T_\beta \phi(x)=\phi T_\beta(x)$ at each simple $\beta$ and $x\in U_q$.
\par

This anti-algebra isomorphism restricts to an isomorphism between $\mfk{v}_{q^{-1}}$ and $\mfk{v}_q$.
\end{lemma}

\begin{proof}
The first claim is immediate from the explicit expressions of the braid group operators \cite[\S 37.1.3]{lusztig93}.  As for $\mfk{v}_q$, the distinguished generators are obtained by applying braid group operators the the $E_\alpha$ and $F_\alpha$ with $l_\alpha>1$.  So one sees that $\phi$ maps the generators for $\mfk{v}_{q^{-1}}$ bijectively to the generators for $\mfk{v}_q$.
\end{proof}

\subsection{Types $B$, $C$ and $F$}

\begin{lemma}
When $G$ is of type $B_n$, $C_n$, or $F_4$, and $l=2$, the algebra $\mfk{v}_q$ is stable under the action of the braid group $B_{\mfk{g}}$.
\end{lemma}

\begin{proof}
The main arguments for types $B_n$ and $C_n$ already appear in the arguments for type $F_4$.  So we only cover the $F_4$ case.
\par

Enumerate the simple roots as $\{\beta_4,\beta_3,\alpha_2,\alpha_1\}$, with the $\beta_i$ long, the $\alpha_i$ short, and consecutive roots neighbors in the Dynkin diagram.  Take $E_i=E_{\alpha_i}$ at $i< 3$ and $E_j=E_{\beta_j}$ at $j\geq 3$.  Consider the subalgebras $\mfk{v}^2_q$, $\mfk{v}^3_q$, $\mfk{v}^4_q$ in $\mfk{v}_q$ defined by
\[
\begin{array}{l}
\mfk{v}^2_q=\text{the subalg gen'd by all $K$'s, $E_1$, $E_2$, $F_1$, and $F_2$}.\vspace{2mm}\\
\mfk{v}^3_q=\text{the subalg gen'd by $\mfk{v}^2_q$ and $E_{32}$, $F_{32}$},\vspace{2mm}\\
\mfk{v}^4_q=\text{the subalg gen'd by $\mfk{v}^3_q$ and $E_{432}$, $F_{432}$}=\mfk{v}_q.
\end{array}
\]
Take also $\mfk{v}^5_q=\mfk{v}_q$.  It suffices to show that $T_\nu(\mfk{v}^j_q)\subseteq \mfk{v}^{j+1}_q$ for each simple $\nu$ and all $j\leq 4$.
\par

We let $T_i=T_{\alpha_i}$, or $T_{\beta_i}$ when appropriate.  In the case $j=2$, $\mfk{v}^2_q$ is immediately seen to be stable under the actions of $T_1$ and $T_2$, and is invariant under $T_4$.  So we need only check $T_3$, and since the $K$'s are stable under the $T_i$ we need only check that $T_3$ applied to the $E$'s and $F$'s lie in $\mfk{v}^3_q$.  The vectors $E_1$ and $F_1$ are $T_3$ invariant, and
\[
T_3(E_2)=E_{32},\ T_3(F_2)=F_{32}\ \in\ \mfk{v}^3_{q}.
\]
We're done.
\par

For the $T_i$ applied to $\mfk{v}^3_{q}$, we need only check the values of these operations on the new vectors $E_{32}$ and $F_{32}$.  One sees directly from the definitions, or the braid group relations, that $T_iT_j(E_r)=T_jT_i(E_r)$ and $T_iT_j(F_r)=T_jT_i(F_r)$ whenever $|i-j|>1$.  So we have
\[
T_1(E_{32})=T_1T_3(E_2)=T_3T_1(E_2)\in T_3(\mfk{v}^2_q)\subseteq \mfk{v}^3_q\subseteq \mfk{v}^4_{q},
\]
and similarly find $T_1(F_{32})\in \mfk{v}^4_{q}$.  Via the rank $2$ calculation \cite[\S 39.2.3 (b)]{lusztig93}, and the fact that $T_3^{-1}=T_3$ by Lemma \ref{lem:T_inv}, we have also
\[
T_2(E_{32})=T_3T_3T_2T_3(E_2)=T_3(E_2)=E_{32}\ \ \text{and similarly}\ \ T_2(F_{32})=F_{32}\in \mfk{v}_q^3.
\]
Here we've used \cite[\S 39.2.3]{lusztig93} specifically for the calculation $T_3T_2T_3(E_2)=E_2$.
\par

For $T_3$, we have $T_3^2=id$ by Lemma \ref{lem:T_inv} so that
\[
T_3(E_{32})=T_3^2(E_2)=E_2\ \ \text{and}\ \ T_3(F_{32})=T_3^2(F_2)=F_2.
\]
For $T_4$ we have simply $T_4(E_{32})=E_{432}$ and $T_4(F_{32})=F_{432}$.  These elements are all in $\mfk{v}^4_{q}$, so that $T_i(\mfk{v}^3_q)\subseteq \mfk{v}^4_{q}$ at all $i$.
\par

For $\mfk{v}^4_q$ we need to check
\[
T_i(E_{432}),\ T_i(F_{432})\ \in \mfk{v}^4_q\ \ \text{for all }i.
\]
For $i=1$ or $2$ we have
\[
T_i(E_{432})=T_iT_4(E_{32})=T_4T_i(E_{32})\in T_4(\mfk{v}^3_q)\subseteq \mfk{v}^4_{q}.
\]
For $i=4$ we have $T_4=T_4^{-1}$ by Lemma \ref{lem:T_inv} so that $T_4(E_{432})=T_4^2(E_{32})=E_{32}\in \mfk{v}^4_q$.  For $i=3$ we apply the braid group relations \cite[Theorem 39.4.3]{lusztig93} to get
\[
T_3(E_{432})=T_3T_4T_3(E_2)=T_4T_3T_4(E_2)=T_4(E_{32})=E_{432}.
\]
The results for $F_{432}$ are obtained similarly.
\end{proof}

\subsection{Type $G_2$}

\begin{lemma}
For $G$ of type $G_2$, and $l=3$, the algebra $\mfk{v}_q$ is stable under the braid group action.
\end{lemma}

\begin{proof}
Write $\Delta=\{\beta,\alpha\}$ with $\alpha$ short and $\beta$ long.  The algebra $\mfk{v}_q$ is generated by $E_\alpha$, $E_{\beta+\alpha}=T_\beta(E_\alpha)$, and the corresponding $F$'s.  We have $T_\beta(E_\alpha)=E_{\beta+\alpha}\in \mfk{v}_q$ already and
\[
T_\beta(E_{\beta+\alpha})=T_\beta^2(E_\alpha)=E_\alpha\ \in\ \mfk{v}_q
\]
by Lemma \ref{lem:T_inv}, since $q_\beta=-1$.  For $\alpha$ we have $T_\alpha(E_\alpha)=-F_\alpha\msf{K}_\alpha$ and
\[
T_\alpha(E_{\beta+\alpha})=E_\alpha^{(2)}E_\beta-q_\alpha^{-2}E_\alpha E_\beta E_\alpha + E_\beta E_\alpha^{(2)}\ \in\ \mfk{v}_q,
\]
by \cite[\S 39.2.2]{lusztig93}.  This shows $T_\nu(\mfk{v}_q^+)\subseteq \mfk{v}_q$ at all simple $\nu$ and, via symmetry, we obtain $T_\nu(\mfk{v}_q^-)\subseteq \mfk{v}_q$ as well.
\end{proof}

\begin{lemma}
For $G$ of type $G_2$, and $l=2$, the algebra $\mfk{v}_q$ is stable under the braid group action.
\end{lemma}

\begin{proof}
We recall
\[
E_{2\alpha+\beta}=E_\alpha^{(2)}E_\beta+E_\beta E_\alpha^{(2)}+E_\alpha E_\beta E_\alpha.
\]
For the simple root vectors we have
\[
T_\alpha(E_\alpha)=-F_\alpha\msf{K}_\alpha,\ \ T_\beta(E_\alpha)=E_{\beta+\alpha}=E_\beta E_\alpha-q_{\alpha}E_\alpha E_\beta,\ \ T_\beta(E_\beta)=-F_\beta \msf{K}_\beta.
\]
We have $[3]_{q_\alpha}=-1$ so that $E_\alpha^{(3)}=-E_\alpha^{(2)}E_\alpha$, which gives
\[
\begin{array}{l}
T_\alpha(E_\beta)=E_\alpha^{(3)}E_\beta+q_\alpha E_\alpha^{(2)}E_\beta E_\alpha - E_\alpha E_\beta E_\alpha^{(2)}-q_\alpha E_\beta E_\alpha^{(3)}\vspace{2mm}\\
\hspace{1cm}=-E_\alpha E_{2\alpha+\beta}+q_\alpha E_{2\alpha+\beta}E_\alpha.
\end{array}
\]
We are only left to check the values of $T_\alpha$ and $T_\beta$ on $E_{2\alpha+\beta}$.
\par

We have
\[
E_\alpha=T_\beta T_\alpha T_\beta T_\alpha T_\beta(E_\alpha)=T_\beta T_\alpha(E_{2\alpha+\beta})
\]
by the rank $2$ calculation \cite[\S 39.2.2 (b)]{lusztig93} and, using the explicit formula for the inverse $T_\beta^{-1}$ \cite[\S 37.1]{lusztig93}, we the obtain
\[
T_\alpha(E_{2\alpha+\beta})=T_\beta^{-1}(E_\alpha)=-q_\alpha E_\beta E_\alpha+E_\alpha E_\beta.
\]
For $T_\beta(E_{2\alpha+\beta})$ we have
\[
T_\beta(E_\alpha^{(2)})=E_\beta^{(2)}E_\alpha^{(2)}-q_\alpha E_\beta E_\alpha^{(2)}E_\beta- E_\alpha^{(2)}E_\beta^{(2)}
\]
\cite[\S 37.1]{lusztig93} and $T_\beta(E_\beta)=\msf{K}_\beta F_\beta$ so that
\[
\begin{array}{l}
T_\beta(E_{2\alpha+\beta}-E_\alpha E_\beta E_\alpha)\vspace{2mm}\\
=\msf{K}_\beta (F_\beta T_\beta(E_\alpha^{(2)})-T_\beta(E_\alpha^{(2)}) F_\beta)\vspace{2mm}\\
=-\msf{K}_\beta ([E_\beta^{(2)},F_\beta]E_\alpha^{(2)}- q_\alpha \binom{\msf{K}_\beta;0}{1}E_\alpha^{(2)}E_\beta- q_\alpha E_\beta E_\alpha^{(2)}\binom{\msf{K}_\beta;0}{1} -E_\alpha^{(2)}[E_\beta^{(2)},F_\beta])\vspace{2mm}\\
=-\msf{K}_\beta([E_\beta^{(2)},F_\beta]E_\alpha^{(2)}-E_\alpha^{(2)}[E_\beta^{(2)},F_\beta])+ q_\alpha \msf{K}_\beta \binom{\msf{K}_\beta;0}{1}(E_\alpha^{(2)}E_\beta+E_\beta E_\alpha^{(2)}).
\end{array}
\]
Now,
\[
[E_\beta^{(2)},F_\beta]=\binom{\msf{K}_\beta;-1}{1}E_\beta
\]
so that the expression reduces to
\[
-q_\alpha \msf{K}_\beta(\binom{\msf{K}_\beta;0}{t}-q_\alpha^{-3} \binom{\msf{K}_\beta;1}{t})(E_\alpha^{(2)}E_\beta+E_\beta E_\alpha^{(2)})
\]
\[
=E_\alpha^{(2)}E_\beta+E_\beta E_\alpha^{(2)}=E_{2\alpha+\beta}-E_\alpha E_\beta E_\alpha
\]
by \cite[\S 6.4 (b4)]{lusztig90II}.  Hence $T_\beta(E_{2\alpha+\beta})=E_{2\alpha+\beta}-E_\alpha E_\beta E_\alpha-q_\alpha \msf{K}_\beta E_{\beta+\alpha}F_\beta E_{\beta+\alpha}$.
\par

The above arguments establish inclusions $T(\mfk{v}_q^+)\subseteq \mfk{v}_q$ for all braid group operators $T$, and by symmetry we have $T(\mfk{v}_q^-)\subseteq \mfk{v}_q$ as well.
\end{proof}

\section{Normality below the lacing number}
\label{sect:norm_app}

We cover the details for the proof of Lemma \ref{lem:norm1} at low order parameters.

\subsection{Completed proof of \ref{lem:norm1}}
\label{sect:norm1}

As argued in the early parts of the proof, we may assume $G$ is almost-simple and that $l>1$.  We first cover the case where $G$ is not of type $G_2$ at $l=2$.

\begin{proof}[Proof away from $G_2$ at $l=2$]
We have that $l$ is equal to the lacing number for $G$, and suppose additionally that the scalar parameter for $G$ is of order $>3$ in type $G_2$.  Equivalently, we assume $q_{\beta}=-1$ whenever $\beta$ is a long simple root.  This implies $K_\beta=K_{\beta}^{-1}$.
\par

It is shown in \cite[Lemma 6.6]{lentner16} that $\mfk{v}_q^+$ remains a normal (braided) Hopf subalgebra in $U_q^+$ when $l$ is equal to the lacing number.  By considering the Hopf isomorphism $U_{q^{-1}}^{\geq 0}\to U_q^{\leq 0}$ which exchanges $E$'s and $F$'s, we also have that $\mfk{v}_q^-$ is normal in $U_q^-$.  Such Hopf normality implies
\begin{equation}\label{eq:2044}
E_\beta^{(l_\beta)}E_\alpha+(-1)^{l_\beta}q(\beta,\alpha)^{l_\beta}E_\alpha E_\beta^{(l_\beta)}\in \mfk{m}^+\ \ \Rightarrow\ [E_\beta^{(l_\beta)},E_\alpha]\ \in\ \hat{U}_q\mfk{m}^+,
\end{equation}
where $\mfk{m}^+$ is the subspace of positive degree elements in $\mfk{v}^+_q$.  But now one sees, by scaling by grouplikes, that $\hat{U}_q\mfk{m}^+\subseteq \hat{U}_q\mfk{m}$.  So we locate the commutator \eqref{eq:2044} in $\hat{U}_q\mfk{m}$.  We similarly find that $[F_\beta^{(l_\beta)},F_\alpha]\in \hat{U}_q\mfk{m}$ at all $\alpha\in \Delta_l$, $\beta\in \Delta$.  So we need only check the commutators of the $E$'s with the $F$'s in \eqref{eq:1045}.
\par

In considering these commutators we have three cases to check; when $\beta$ is short and has no long neighbors, when $\beta$ is short but has a long neighbor, and when $\beta$ is long.  In the last case $l_\beta=1$.  Let us list the simple roots as
\[
\Delta=\{\beta_n,\dots,\beta_{r+1},\alpha_r,\dots,\alpha_1\}
\]
\[
\Delta_l=\{\beta_n+\cdots+\beta_{r+1}+\alpha_r,\dots,\beta_{r+1}+\alpha_r,\alpha_r,\dots,\alpha_1\}
\]
with neighbors appearing in consecutive order, all the $\alpha_i$ short, and all $\beta_j$ long.  Take
\[
E_i=E_{\alpha_i}\ \text{or}\ E_{\beta_i}\ \text{where appropriate, }E_{r,r}=E_r\text{ and }E_{r+j,r}=E_{\beta_{r+j}+\cdots+\beta_{r+1}+\alpha_r}.
\]
Adopt a similar notation for the $F$'s.

Recall that $q_\beta=-1$ whenever $\beta$ is long.  When $\beta$ is short and has no long neighbors $[F_\beta^{(l_\beta)},E_{r+j,r}]=[E_\beta^{(l_\beta)},F_{r+j,r}]=0$ for all positive $j$.  So the computation reduces to a type $A_n$ calculation at $l=2$, which we have already addressed in Section \ref{sect:normal.1}.
\par

When $\beta$ is short and has a long neighbor we have $\beta=\alpha_r$ necessarily.  The generic commutator relation \cite[\S 6.5 (a2)]{lusztig90II} shows that $[E_j,F_\beta^{(l_\beta)}]$ lies in $\hat{U}_q\mfk{m}$ for all $j$.  For $E_{r+1,r}$ we have
\[
[E_{r+1,r},F_\beta^{(l_\beta)}]=E_{r+1}[E_r,F_{\beta}^{(l_\beta)}]+[E_r,F_{\beta}^{(l_\beta)}]E_{r+1}.
\]
Now
\[
[E_r,F_{\beta}^{(l_\beta)}]=\binom{K_\beta;m}{1}F_\beta^{(l_\beta-1)}
\]
for some integer $m$ \cite[\S 6.5 (a2), (a6)]{lusztig90II}.  Since $F_\beta$ commutes with $E_{r+1}$ in this case we get
\[
[E_{r+1,r},F_\beta^{(l_\beta)}]=(\binom{K_\beta;m'}{1}+\binom{K_\beta;m}{1})E_{r+1}F^{(l_\beta-1)}_\beta.
\]
Since $F_\beta^{(l_\beta-1)}$ is in $\mfk{m}\subseteq \bar{u}_q$ this element is in $\hat{U}_q\mfk{m}$, as desired.  One now sees by induction that
\[
[E_{r+j,r},F_\beta^{(l_\beta)}]\ \in\ \hat{U}_q F_\beta\ \subseteq\ \hat{U}_q\mfk{m}
\]
for all positive $j$, when $\beta=\alpha_r$.

We are left, finally, to deal with the case where $\beta$ is long, so that $l_\beta=1$.  We have $F_\beta^{(l_\beta)}=F_{r+j}$ for some positive $j$.  Now, when $j>j'\geq 0$ we have immediately $[F_{r+j},E_{r+j',r}]=0$.  When $j'=j$ we have
\[
[F_{r+j},E_{r+j,r}]=-[F_{r+j},E_{r+j}E_{r+j-1,r}+E_{r+j-1,r}E_{r+j}]
\]
\begin{equation}\label{eq:1112}
=-\binom{K_{r+j};0}{1}E_{r+j-1,r}-E_{r+j-1,r}\binom{K_{r+j};0}{1}.
\end{equation}
Since $q_{\beta}=-1$, we have the commutativity relations
\[
E_{r+j-1}\binom{K_{r+j};0}{1}=\binom{K_{r+j};1}{1}E_{r+j-1}=-\binom{K_{r+j};0}{1}E_{r+j-1}+K_\beta E_{r+j-1}
\]
\cite[\S 6.4 and 6.5 (b5), (a5)]{lusztig90II} which reduces \eqref{eq:1112} to
\[
[F_{r+j},E_{r+j,r}]=-K_{r+j} E_{r+j-1,r}\ \in\ \mfk{m}.
\]
For $r+j+1$ the above expression gives also
\[
\begin{array}{l}
[F_{r+j},E_{r+j+1,r}]=-E_{r+j+1}K_{r+j}E_{r+j-1,r}-K_{r+j}E_{r+j-1,r}E_{r+j+1}\vspace{2mm}\\
\ \ =-(q_{\beta_{r+j}}^{-1}+1)K_{r+j}E_{r+j-1,r}E_{r+j+1}=0,
\end{array}
\]
since $q_{\beta_{r+j}}=-1$.  By induction one now finds $[F_{r+j},E_{r+j+k,r}]=0$ at all $k\geq 1$.  We have now checked all cases and obtain the desired inclusion
\[
[F_\beta^{(l_\beta)},E_\alpha]\ \in\ \hat{U}_q\mfk{m}
\]
for all $\beta\in \Delta$ and $\alpha\in \Delta_l$.  The cases with the $E$'s and $F$'s swapped follow by symmetry.
\par

We have now established the desired inclusions \eqref{eq:1045} whenever $G$ is not of type $G_2$ at a $3$-rd root of unity.  The final case of $G_2$ at a $3$-rd root of $1$ is established via completely similar arguments, though some signs are changed.
\end{proof}

We now cover the special case of $G_2$ at a $4$-th root of $1$.

\begin{proof}[Proof for $G_2$ at $l=2$]
As in the analysis at \eqref{eq:2044}, we can use Lentner's normality result \cite[Lemma 6.6]{lentner16} to find that
\[
E_\alpha E_\beta^{(l_\nu)},\ F_\alpha F_\beta^{(l_\nu)}\ \in \hat{U}_q\mfk{m}
\]
whenever $\alpha\in \Delta_l$ and $\beta\in\Delta$.  So we need only check the products $E_\alpha F_\beta^{(l_\nu)}$, $F_\alpha E_\beta^{(l_\nu)}$.
\par

In this special setting $l_\beta=2$ at all simple $\beta$.  We have
\[
[E_\beta^{(l_\beta)},F_{\alpha}],\ [F_\beta^{(l_\beta)},E_{\alpha}]\ \in\ \hat{U}_q\mfk{m} 
\]
whenever $\alpha,\beta\in \Delta$ by the generic commutativity relation \cite[\S 6.5 (a2)]{lusztig93}.  For the unique non-simple root $2\alpha+\beta$ in $\Delta_l$ we have
\begin{equation}\label{eq:1167}
E_{2\alpha+\beta}=
E_\alpha^{(2)} E_\beta+E_\beta E_\alpha^{(2)}+E_\alpha E_\beta E_\alpha
\end{equation}
and compute
\[
[F_\beta^{(2)},E_\alpha^{(2)}]=[F_\beta^{(2)},E_\alpha]=0,\ \ [F^{(2)}_\beta,E_\beta]=\binom{K_\beta;-1}{1}F_\beta,\ \ [F_\alpha^{(2)},E_\beta]=0
\]
\[
[F_\alpha^{(2)},E_\alpha]=\binom{K_\alpha;-1}{1}F_\beta,\ \ \text{and}\ \ [F_\alpha^{(2)},E_\alpha^{(2)}]=\binom{K_\alpha;m}{1}F_\alpha E_\alpha-\binom{K_\alpha;0}{2},
\]
where $m$ is some integer \cite[\S 6.5]{lusztig90II}.  Hence applying the derivation $[F_\beta^{(2)},-]$ to the expression \eqref{eq:1167} yields an elements
\[
[F_\beta^{(2)},E_{2\alpha+\beta}]\ \in\ \hat{U}_q F_\beta\ \subseteq\ \hat{U}_q\mfk{m}.
\]
This implies $E_{2\alpha+\beta}F_\beta^{(2)}\in \hat{U}_q\mfk{m}$.

As for the derivation $[F_\alpha^{(2)},-]$, we have
\[
E_\beta\binom{K_\alpha;0}{2}=\binom{K_\alpha;2}{2}E_\beta
\]
\cite[\S 6.5]{lusztig90II} so that we find again
\[
[F_\alpha^{(2)},E_{2\alpha+\beta}]\ \in\ \hat{U}_q\mfk{m}\ \Rightarrow\ E_{2\alpha+\beta}F_\alpha^{(2)}\in \hat{U}_q\mfk{m}.
\]
We employ symmetry to obtain $F_{2\alpha+\beta}E_\beta^{(2)},\ F_{2\alpha+\beta}E_\alpha^{(2)}\in \hat{U}_q\mfk{m}$ as well.  This completes the proof.
\end{proof}

\end{document}